\documentclass[12pt]{amsart}
\usepackage{amssymb,latexsym}
\usepackage{enumerate}
\usepackage{amscd}
\usepackage{verbatim}

\theoremstyle{plain}
\newtheorem{theorem}{Theorem}
\newtheorem{corollary}[theorem]{Corollary}
\newtheorem{lemma}[theorem]{Lemma}
\newtheorem{proposition}[theorem]{Proposition}

\theoremstyle{definition}

\theoremstyle{remark}
\newtheorem*{remark}{Remark}
\newtheorem*{notation}{Notation}

\newtheorem*{example}{Example}

\numberwithin{theorem}{section}
\numberwithin{definition}{section}

\def\pair#1#2{\langle #1, #2\rangle}

\def\Mof#1{\Cal M(A)}

\newcommand{\N}{{\mathbb N}}
\newcommand{\Z}{{\mathbb Z}}

\def\Unif{{\text{\bf Unif}}}

\newcommand{\Cg}{\operatorname{Cg}}
\newcommand{\Ug}{\operatorname{Ug}}
\newcommand{\Con}{\operatorname{Con}}

\newcommand{\Fg}{\operatorname{Fg}}

\newcommand{\OpSemiUnif}{\operatorname{SemiUnif}}

\newcommand{\nat}{\operatorname{nat}}

\newcommand{\Fil}{\operatorname{Fil}}

\newcommand{\I}{\operatorname{I}}
\newcommand{\OpUnif}{\operatorname{Unif}}

\newcommand{\Rel}{\operatorname{Rel}}

\def\congruence{on\-gru\-ence\discretionary{-}{}{-}}
\def\conM/{c\congruence mod\-u\-lar}
\def\ConM/{C\congruence mod\-u\-lar}
\def\conD/{c\congruence dis\-trib\-u\-tive}
\def\ConD/{C\congruence dis\-trib\-u\-tive}
\def\conP/{c\congruence per\-mut\-a\-ble}
\def\ConP/{C\congruence per\-mut\-a\-ble}
\def\conMity/{\conM/\-i\-ty}
\def\ConMity/{\ConM/\-i\-ty}
\def\conDity/{c\congruence dis\-trib\-u\-tiv\-i\-ty}
\def\ConDity/{C\congruence dis\-trib\-u\-tiv\-i\-ty}
\def\conPity/{c\congruence per\-mut\-a\-bil\-i\-ty}
\def\ConPity/{C\congruence per\-mut\-a\-bil\-i\-ty}
\def\usprv/{un\-der\-ly\-ing-set-pre\-ser\-ving}

\def\ie/{{i.e.}}
\def\Ie/{{I.e.}}
\def\eg/{{e.g.}}
\def\Eg/{{E.g.}}
\def\etc/{{etc.}}

\newdimen\mysubdimen
\newbox\mysubbox

\def\subwhat#1#2#3{{
\setbox\mysubbox=\hbox{#3}
\mysubdimen=\wd\mysubbox
\setbox\mysubbox=\hbox{$#1#2$}
\ifnum\mysubdimen>\wd\mysubbox
\vtop{
\hbox to\mysubdimen{\hfil\box\mysubbox\hfil}
\nointerlineskip
\hbox{#3}}
\else
\mysubdimen=\wd\mysubbox
\vtop{
\box\mysubbox
\nointerlineskip
\hbox to\mysubdimen{#3}}
\fi
}}

\begin{document}

\title{Commutator Theory for Compatible Uniformities}
\author{William H. Rowan}
\address{PO Box 20161 \\
         Oakland, California 94620}
\email{whrowan@member.ams.org}


\date{\today}

\begin{abstract}
We investigate commutator operations on compatible
uniformities. We define a commutator operation for
uniformities in the congruence-modular case which extends
the commutator on congruences, and explore its properties.
\end{abstract}

\maketitle

\section*{Introduction}

The purpose of this paper is to generalize the commutator
of congruences to a commutator of compatible
uniformities. Commutator theory (on congruences) works
best for congruences of algebras in congruence-modular
varieties. The same is true of the commutator of
uniformities described here. In fact, the commutator
of congruences $\alpha$ and
$\beta$ becomes a special case of that of uniformities,
when we view
$\alpha$ and
$\beta$ as the uniformities $\Ug\{\,\alpha\,\}$ and
$\Ug\{\,\beta\,\}$ that they generate, because we have
$\Ug\{\,[\alpha,\beta]\,\}=[\Ug\{\,\alpha\,\},
\Ug\{\,\beta\,\}]$.

We follow the development of Commutator Theory in
\cite{F-McK} fairly closely. The main thesis of
\cite{R02}, where compatible uniformities were first
studied systematically in the context of Universal
Algebra, is that compatible uniformities can be considered
a generalization of congruences. Often, there is
a reasonably direct translation of congruence-theoretic
arguments into uniformity-theoretic ones. Following this
philosophy, we are able to generalize (in
Sections~\ref{S:Term} and
\ref{S:Commutator}) the concept of
$C(\alpha,\beta;\delta)$ ($\alpha$ \emph{centralizes}
$\beta$ \emph{modulo} $\delta$) to compatible
uniformities, and in the congruence-modular case, to
define $[\mathcal U,\mathcal V]$ to be the least
uniformity $\mathcal W$ such that
$C(\mathcal U,\mathcal V;\mathcal W)$.

Another approach to the commutator $[\alpha,\beta]$, for
congruences $\alpha$ and $\beta$, as discussed in
\cite{F-McK}, is to study
congruences of the algebra $A(\alpha)$. The congruence
$\beta$ is pushed out along the homomorphism
$\Delta_\alpha:A\to A(\alpha)$ that sends $a\in A$ to
$\pair aa$, yielding a congruence $\Delta_{\alpha,\beta}$
which gives rise to
$[\alpha,\beta]$. In the case of uniformities, we can
replace $\beta$ by a uniformity
$\mathcal U$, and push it out along $\Delta_\alpha$,
yielding a compatible uniformity $\Delta_{\alpha,\mathcal
U}$ on
$A(\alpha)$ which we then show gives rise to
$[\Ug\{\,\alpha\,\},\mathcal U]$ in the
important special case of algebras having term operations
comprising a group structure. (This includes
many familiar varieties of algebras, such as groups,
rings, and varieties of nonassociative algebras.)
This is done in Section~\ref{S:A(alpha)}.  It is natural
to ask whether the theory can be extended to give an
interpretation of
$[\mathcal U,\mathcal V]$ in terms of
compatible uniformities on some algebra $A(\mathcal U)$.
Unfortunately, the necessary definition of $A(\mathcal
U)$ is not yet available.

In the congruence-modular case, the properties of
the commutator on compatible uniformities described in
Section~\ref{S:Properties} duplicate those of the
commutator of congruences, with some regrettable gaps in
what we have been able to prove. In particular, because
of the difficulties we encounter in working with joins of
uniformities, the additivity of the commutator is not
settled in general, although we prove it is true for some
important special cases.

We also describe a
natural commutator operation for congruential
uniformities, which are uniformities that, as a filter of
relations, have a base of congruences, and are thus of
the form $\Ug F$ for a filter $F$ in $\Con A$. The
commutator $[\Ug F,\Ug F']$ is defined in terms of the
commutators of elements of $F$ and $F'$; see
Section~\ref{S:Congruential}.

Section~\ref{S:Examp} is devoted to miscellaneous matters,
including commutators of congruential uniformities on
commutative rings. We prove that in that case, the two
definitions of the commutator on congruential
uniformities, the one given by the general definition of
Section~\ref{S:Commutator} and the other given by the
formula of Section~\ref{S:Congruential}, coincide. This
appears to be a special property of commutative rings; in
general, we do not know even whether the commutator
operation of Section~\ref{S:Commutator}, applied to
congruential uniformities, always gives a congruential
uniformity. We also discuss in this section
the case of varieties which are
congruence-distributive, where we show that as in the
case of the commutator of congruences, the commutator of
two compatible uniformities is simply their meet.

In the last section, we
discuss the current state of some questions about
compatible uniformities, uniformity lattices, and
commutators of uniformities.

\section*{Preliminaries}

\subsection{Category theory} We follow \cite{MacL} in
terminology and notation.  In particular, $1_a$ will stand
for the identity arrow on an object $a$ in a category
$\mathbf C$.

\subsection{Lattice theory} The reader should be
familiar with lattices. We use $\top$ and
$\bot$ to denote the greatest and least elements of a
lattice, assuming they exist, and
$\wedge$ and
$\vee$ for the meet and join operations.

\subsection{Filters}

If $L$ is a lattice, then a nonempty subset $F\subseteq
L$ is called a \emph{filter} if $y\geq x\in F$ implies
$y\in F$ and $x$, $y\in F$ imply $x\wedge y\in F$.

If $S\subseteq L$ is a nonempty set, then the
\emph{filter generated by $S$}, denoted by $\Fg^L S$
or simply $\Fg S$, is the subset of elements of $L$ that
are greater than or equal to a finite meet of elements of
$S$. An important special case, given $x\in L$, is
$\Fg\{\,x\,\}$, the
\emph{principal filter generated by $x$}, which is just
the set of elements of $L$ greater than or equal to $x$.

Filters are ordered by reverse inclusion and the filters
in a nonempty lattice form a complete lattice. The meet
of a tuple of filters $F_i$ is given
by $\bigwedge_iF_i=\Fg(\bigcup_iF_i)$. The join of the
tuple is the intersection: $\bigvee_iF_i=\bigcap_iF_i$.

If $F$ is a filter, a \emph{base} for $F$ is a subset
$B\subseteq F$ such that $x\in F$ implies $b\leq x$ for
some $b\in B$. If $L$ is a lattice, then a subset
$B\subseteq L$ is a base for a filter of $L$, or
\emph{filter base}, iff given any $x$, $y\in B$, there is
a
$z\in B$ such that
$z\leq x\wedge y$.

\subsection{Universal algebra}

We assume familiarity with universal algebra, as
explained for example in \cite{B-S}.  We prefer to allow
an algebra to have an empty underlying set, however. We
denote the underlying set of an algebra $A$ by $|A|$.

If $R$ is a binary relation on (the underlying
set of) an algebra $A$, then we denote by
$\Cg R$ the smallest congruence $\alpha\in\Con A$ such
that
$R\subseteq\alpha$.

\subsection{Comgruence-permutable, congruence-modular,
and congruence-distributive varieties}

A variety of algebras $\mathbf V$ is
\emph{congruence-permutable} (or, a \emph{Mal'tsev
variety}) if for every
$A\in\mathbf V$, for every $\alpha$, $\beta\in\Con A$,
$\alpha\circ\beta=\beta\circ\alpha$. A variety is
congruence-permutable iff \cite{Mal} there is a
ternary term $p$, called a \emph{Mal'tsev term},
satisfying the identities
$p(x,x,y)=y$ and $p(x,y,y)=x$. For example, the variety
of groups is congruence-permutable because it has the
Mal'tsev term $p(x,y,z)=xy^{-1}z$.

 A
variety of algebras
$\mathbf V$ is
\emph{congruence-modular} if for every $A\in\mathbf V$,
$\Con A$ is a modular lattice. A variety is
congruence-modular iff \cite{Day} there is a
finite sequence $m_0$, $\ldots$, $m_k$ of quaternary
terms, called \emph{Day terms}, satisfying the following
identities:
\begin{enumerate}
\item[(D1)] For all $i$, $m_i(x,y,y,x)=x$
\item[(D2)] $m_0(x,y,z,w)=x$
\item[(D3)] $m_k(x,y,z,w)=w$
\item[(D4)] for even $i<k$,
$m_i(x,x,y,y)=m_{i+1}(x,x,y,y)$
\item[(D5)] for odd $i<k$,
$m_i(x,y,y,z)=m_{i+1}(x,y,y,z)$.
\end{enumerate}

A congruence-permutable variety, with Mal'tsev term $p$,
is necessarily congruence-modular, with Day terms
$m_0(x,y,z,w)=x$,
$m_1(x,y,z,w)=p(x,p(x,y,z),w)$, and $m_2(x,y,z,w)=w$.

Similarly, a variety of algebras
$\mathbf V$ is
\emph{congruence-distributive} if for every $A\in\mathbf
V$, $\Con A$ is distributive. A variety is
congruence-distributive iff
\cite{Jonsson} there is a finite sequence $d_0$,
$\ldots$, $d_k$ of ternary terms, satisfying the
following identities:
\begin{enumerate}
\item[(J1)] $d_0(x,y,z)=x$
\item[(J2)] $d_n(x,y,z)=z$
\item[(J3)] for all $i$, $d_i(x,y,x)=x$
\item[(J4)] for even $i<k$, $d_i(x,x,y)=d_{i+1}(x,x,y)$
\item[(J5)] for odd $i<k$, $d_i(x,y,y)=d_{i+1}(x,y,y)$.
\end{enumerate}

\subsection{Commutator theory}

The \emph{commutator} is a binary operation on
congruences in the congruence lattice $\Con A$ of an
algebra $A$ in a congruence-modular variety, and which
is sometimes defined for more general varieties.

If
$A$ is an algebra in a congruence-modular variety, and
$\alpha$,
$\beta\in\Con A$, we can define the
commutator $[\alpha,\beta]$ as the least
$\delta\in\Con A$ such that $\alpha$ \emph{centralizes}
$\beta$ \emph{modulo}
$\delta$, or in other words, such that $\delta$ satisfies
the \emph{$\alpha,\beta$-term condition}. See the first
part of Section~\ref{S:Term} for detailed definitions. 
Whereas Section~\ref{S:Term} gives these definitions in
the general case, and doesn't really define
$[\alpha,\beta]$ or its generalization to uniformities,
Section~\ref{S:Commutator} makes the assumption of
congruence-modularity and proves the simplifications
that make the definitions of
$[\alpha,\beta]$ and $[\mathcal U,\mathcal V]$ (for
$\mathcal U$, $\mathcal V$ compatible uniformities) so
reasonable.

The commutator is so named because it generalizes the
notion of the commutator of normal subgroups of a group.
(Of course, the variety of groups is
congruence-modular.)  As a further example, the variety
of commutative rings is congruence-modular, and for a
commutative ring
$A$, the commutator is simply the product of ideals. That
is, if
$I_\alpha$ denotes the ideal corresponding to
$\alpha\in\Con A$, then we have
$I_{[\alpha,\beta]}=I_\alpha I_\beta$.

\subsection{Uniform universal algebra}
Universal algebra over the base category $\Unif$ of
uniform spaces, as opposed to the category of sets, was
first studied systematically in
\cite{R02}. This paper develops commutator theory as a
part of that subject. It will be best if the reader has
access to \cite{R02} while reading this paper, but we
will also devote most of the next two sections to
summarizing some of the basic definitions and results
that we need.

\section{Uniformities}
\label{S:Unif}

We denote the set of binary relations on a set $S$ by
$\Rel S$. $\Rel S$, ordered by inclusion, is a complete
lattice.

If $R\in\Rel S$, then by $R^{-1}$ we mean the relation
$\{\,\pair xy\mid\pair yx\in R\,\}$, and by $R^n$, for
$n>0$, we mean the $n$-fold relational product of $n$
copies of $R$.

Consider the following five conditions on a set
$\mathcal U\subseteq\Rel S$:
\begin{enumerate}
\item[(U1)] if $U\in\mathcal U$, and $U\subseteq V$, then
$V\in\mathcal U$
\item[(U2)] if $U$, $V\in\mathcal U$, then $U\cap
V\in\mathcal U$
\item[(U3)] if $U\in\mathcal U$, then
$\Delta_S\overset{\text{def}}=\{\,\pair ss\mid s\in
S\,\}\subseteq U$
\item[(U4)] if $U\in\mathcal U$, then $U^{-1}\in\mathcal
U$
\item[(U5)] if $U\in\mathcal U$, then $V^2\subseteq
U$ for some $V\in\mathcal U$.
\end{enumerate}
Then we say that $\mathcal U$ is a \emph{semiuniformity}
if $\mathcal U$ satisfies conditions (U1) through (U4),
and a \emph{uniformity} if it satisfies (U1) through (U5).
Note that conditions (U1) and (U2) simply state that
$\mathcal U$ is a filter of binary relations

\begin{proposition} A filter $\mathcal U$ satisfies (U3) 
iff
$\Fg^{\Rel S}\{\,\Delta_S\,\}\leq\mathcal U$.
\end{proposition}

\begin{proposition} A filter $\mathcal U$ satisfies (U4)
iff
$\mathcal U=\mathcal U^{-1}$, where $\mathcal
U^{-1}=\{\,U^{-1}\mid U\in\mathcal U\,\}$.
\end{proposition}

\begin{proposition} If $\mathcal U$, $\mathcal V$ are
filters in $\Rel S$, then $\{\,U\circ V\mid U\in\mathcal
U,\,V\in\mathcal V\,\}$ is a base for a filter $\mathcal
U\circ\mathcal V$ in $\Rel S$.
\end{proposition}

Note that for filters $\mathcal U$ and $\mathcal V$ of
reflexive relations, $\mathcal U\cap\mathcal
V\leq\mathcal U\circ\mathcal V$. Thus, $\mathcal
U\leq\mathcal U\circ\mathcal U$ if $\mathcal U$ satisfies
(U3).

\begin{proposition} A filter $\mathcal U$ satisfies (U5)
iff
$\mathcal U\circ\mathcal U\leq\mathcal U$.
\end{proposition}

\begin{notation} If $U\in\mathcal U$, where $\mathcal U$
satisfies (U5), then by induction we can show that there
is a $V\in\mathcal U$ such that $V^n\subseteq U$. We denote such a
$V$ by
$^nU$. This notation must be used with care, particularly
in relation to quantifiers; we do not mean that
$^nU$ is a function of
$U$; it is simply a shorthand for the statement that
there exists such a $V$ and that we will denote it by
$^nU$.
\end{notation}

\subsection{The lattice operations}

We denote the set of uniformities on a set $S$ by
$\OpUnif S$, and the set of semiuniformities by
$\OpSemiUnif S$. We order
these sets by reverse inclusion, i.e., the ordering
inherited from
$\Fil\Rel S$.

The meet of an arbitrary tuple of uniformities on $S$, in
the lattice
$\Fil\Rel S$, is a uniformity. Thus, $\OpUnif S$ admits
arbitrary meets, and is a complete lattice. The same is
true for $\OpSemiUnif S$. The join of a tuple of
semiuniformities is simply the intersection, and
$\OpSemiUnif S$ is a distributive lattice. The theory of
joins of uniformities is more difficult.

\subsection{Permutability} Permutability of congruences
is an important condition in Universal Algebra, and the
condition can be generalized to uniformities. Note that
this subject was first discussed in \cite{R02}, but the
discussion there is not entirely correct; in particular,
Theorem 6.1 is wrong.

If $\mathcal U$ and $\mathcal V$ are uniformities on a
set $S$, then we say that $\mathcal U$ and $\mathcal V$
\emph{permute} if $\mathcal U\circ\mathcal V=\mathcal
V\circ\mathcal U$, and that $\mathcal U$ and $\mathcal V$
\emph{semipermute} if either $\mathcal U\circ\mathcal
V\leq\mathcal V\circ\mathcal U$, or $\mathcal
V\circ\mathcal U\leq\mathcal U\circ\mathcal V$.

As we mentioned, the join operation in the lattice of
uniformities can be difficult to deal with in the general
case, but the case where $\mathcal U$ and $\mathcal V$
semipermute is an easy and important special case; the
following theorem is a revised and corrected version of
\cite[Theorem~6.1]{R02}:

\begin{theorem} Let $\mathcal U$, $\mathcal V\in\OpUnif
S$. Then $\mathcal U\vee\mathcal V=\mathcal
U\circ\mathcal V$ iff $\mathcal
V\circ\mathcal U\leq\mathcal U\circ\mathcal V$.
\end{theorem}

\begin{proof} ($\impliedby$): It is trivial that $\mathcal
U\circ\mathcal V$ satisfies (U1), (U2), and (U3). If
$U\in\mathcal U$ and $V\in\mathcal V$, then since
$\mathcal V\circ\mathcal U\leq\mathcal U\circ\mathcal
V$, there are $\bar U\in\mathcal U$ and $\bar
V\in\mathcal V$ such that $\bar V\circ\bar U\subseteq
U\circ V$. But, $\bar V\circ\bar U=((\bar
U^{-1})\circ(\bar V^{-1}))^{-1}$. Thus, $\mathcal
U\circ\mathcal V$ also satisfies (U4).
Finally, $(\mathcal U\circ\mathcal V)\circ(\mathcal
U\circ\mathcal V)\leq(\mathcal U\circ\mathcal
U)\circ(\mathcal V\circ\mathcal V)\leq\mathcal
U\circ\mathcal V$, verifying (U5). Thus, $\mathcal
U\circ\mathcal V$ is a uniformity. Since $\mathcal
U\leq\mathcal U\circ\mathcal V$ and $\mathcal
V\leq\mathcal U\circ\mathcal V$, we have $\mathcal
U\vee\mathcal V\leq\mathcal U\circ\mathcal V$.
However, $\mathcal U\circ\mathcal V\leq\mathcal U\vee
\mathcal V=\mathcal U\circ\mathcal V$ by (U5).

 ($\implies$): If $\mathcal U\vee\mathcal V=\mathcal
U\circ\mathcal V$, then $\mathcal V\circ\mathcal
U\leq\mathcal U\vee\mathcal V$ by (U5).

\end{proof}

The results in \cite{R02} that
use permutability as a hypothesis, except for Theorem 6.1,   are correct, and
remain true if the hypothesis is weakened to
semipermutability.

\section{Compatible Uniformities}
\label{S:Compat}

\subsection{Compatibility} If $R$ is a relation on an
algebra $A$, we say that $R$ is \emph{compatible} (with
the operations of $A$) if $\mathbf a\mathrel R\mathbf
a'$ implies $\omega^A(\mathbf a)\mathrel R
\omega^A(\mathbf a')$ for each operation symbol $\omega$,
where
$\mathbf a\mathrel R\mathbf{a'}$ means that $a_i\wedge
a'_i$ for all $i$.

We say
that a filter
$\mathcal U$ of reflexive relations on an algebra $A$ is
\emph{compatible} if for each $U\in\mathcal U$, and
each basic operation symbol $\omega$,
there is a $\bar U\in\mathcal U$ such
that $\omega(\bar U)=\{\,\pair{\omega(\mathrm
x)}{\omega(\mathbf y)}\mid x_i\mathrel{\bar
U}y_i\text{ for all }i\,\}\subseteq U$. In this case, for
any term
$t$, given $U\in\mathcal U$, there is a $\bar
U\in\mathcal U$ such that $t(\bar U)\subseteq U$.

We say that $\mathcal
U$ is \emph{singly compatible} if for each $n$-ary term
$t$ for $n\geq 1$, given $U\in\mathcal U$, there is a
$\bar U\in\mathcal U$ such that $t(\bar U,\mathbf
a)\subseteq U$ for every $\mathbf a\in A^{n-1}$.

\begin{lemma}\label{T:WCompatU} If $\mathcal U$ is a
uniformity on $A$, which is singly compatible, then
$\mathcal U$ is compatible.
\end{lemma}

\begin{proof} This follows easily from (U5).
\end{proof}

As a result of the Lemma, single compatibility will be of
interest to us only for semiuniformities.

If $A$ is an algebra, we denote by $\OpSemiUnif A$ the
set of compatible semiuniformities on $A$, and by
$\OpUnif A$ the set of compatible uniformities.

\begin{remark} Since a set $S$ can be seen as an algebra
with no operations, the theory of $\OpUnif S$ is subsumed
by the theory of $\OpUnif A$ for an algebra $A$. Parts of
this section are therefore pertinent the study of
$\OpUnif S$ where $S$ is just a set.
\end{remark}

If filters $\mathcal U_i$ are compatible, so is
$\bigwedge_i\mathcal U_i$. It follows that the meet of
an arbitrary tuple of compatible uniformities or
semiuniformities is also compatible, so the sets of
compatible uniformities and compatible semiuniformities
are complete lattices.  Similarly, if $\mathcal U_i$ is a
tuple of singly compatible semiuniformities, then
$\bigwedge_i\mathcal U_i$ is a singly compatible
semiuniformity.

\subsection{$\Ug\mathcal R$}

If $\mathcal R$ is a filter of relations on an algebra
$A$, then $\Ug\mathcal R$ will denote the smallest
compatible uniformity $\mathcal U$ such that $\mathcal
R\leq\mathcal U$. If we mean instead the smallest
not-necessarily compatible uniformity, we will write
$\Ug^{|A|}\mathcal R$.

\subsection{Joins}

\begin{proposition} If $\mathcal U_i$, $i\in I$ are
elements of $\Unif A$, then $\bigvee_i\mathcal
U_i=\Ug(\bigcap_i\mathcal U_i)$.
\end{proposition}

\begin{theorem} Let $A$ be an algebra.
We have
\begin{enumerate}
\item[(1)]
The join (in the lattice $\OpUnif|A|$) of a tuple of
compatible uniformities is compatible.
\item[(2)] The join (in the lattice $\OpSemiUnif|A|$) of a
tuple of singly compatible semiunifomities is singly
compatible.
\end{enumerate}
\end{theorem}

\begin{proof} (1): See \cite[Theorem~5.3]{R02}.

(2):
If $\mathcal U_i$ are a tuple of singly compatible
semiuniformities on $A$, then their join (in the lattice
of semiuniformities on $|A|$) is $\mathcal
U=\bigcap_i\mathcal U_i$. Let $t$ be an $n$-ary term
operation. If $U\in\mathcal U$, then for each $i$, there
is $U_i\in\mathcal U_i$ with $U_i\subseteq U$. Since
each $\mathcal U_i$ is singly compatible,
 $t(\bar
U_i,\Delta_A,\ldots,\Delta_A)\subseteq U$ for some
$\bar U_i\in\mathcal U_i$. Then $\bigcup_iU_i\in\mathcal
U$, and $t(\bigcup_i\bar
U_i,\Delta_A,\ldots,\Delta_A)\subseteq U$.
\end{proof}

\subsection{Examples of compatible uniformities}

An important special case of a compatible uniformity is
given by choosing
$\mathcal R=\{\,\alpha\,\}$ where $\alpha\in\Con A$. Then
$\Ug\mathcal R=\Fg^{\Rel A}\{\,\alpha\,\}$.

More generally, we can consider $\Ug\{\,\rho\,\}$ where
$\rho\in\Rel A$. However, we have

\begin{theorem} $\Ug\{\,\rho\,\}=\Ug\{\,\Cg\rho\,\}$, and
is compatible if $\rho$ is.
\end{theorem}

\begin{proof}
It suffices to show that
$\Ug\{\,\Cg\rho\,\}\subseteq\Ug\rho$, or in other words
that if $U\in\Ug\{\,\rho\,\}$, then
$\Cg\{\,\rho\,\}\subseteq U$.

Let $U\in\Ug\{\,\rho\,\}$. We have $\rho\subseteq U$, so
$\rho\cup\Delta\subseteq U$ by (U3). Thus, we can reduce
to the case where $\rho$ is reflexive by replacing $\rho$
with $\rho\cup\Delta$.

If $U\in\Ug\{\,\rho\,\}$, then $\rho\subseteq U^{-1}$, so
$\rho^{-1}\subseteq U$. Thus, we can further reduce to
the case where $\rho$ is symmetric, by replacing $\rho$
by $\rho\cup\rho^{-1}$.

Finally, if $U\in\Ug\{\,\rho\,\}$, then $\rho\subseteq
{^nU}$ for all $n\in\mathbb N$, which implies
$\rho\subseteq U$, showing that
$\Cg\rho=\bigcup_n\rho^n\subseteq U$.

As regards compatibility, it is easy to prove that if
$\rho$ is compatible, then so is $\Cg\rho$. It is obvious
that $\Ug\{\,\alpha\,\}$ is compatible if $\alpha$ is a
congruence.
\end{proof}

Another important special case is $\mathcal R=F$, where
$F$ is a filter in $\Con A$. In this case, $\Ug\mathcal
R=\Fg^{\Rel A}F$. A uniformity of this form called a
\emph{congruential uniformity}.

\subsection{Uniformities and congruences}

If $\mathcal U\in\OpUnif A$, then $\bigcap\mathcal
U\in\Con A$. We may consider this as
a mapping from $\OpUnif A$ to $\OpUnif A$, where we map
$\mathcal U$ to $\Ug\{\,\bigcap\mathcal U\,\}$; more
generally, we can map $\mathcal U$ to the filter of
$\kappa$-fold intersections of relations in $\mathcal U$,
for $\kappa$ some given infinite cardinal. The result
will be a compatible uniformity $\mathcal V$ such that
$\mathcal V$ admits $\kappa$-fold intersections of its
elements. We say that $\mathcal V$ \emph{satisfies the
$\kappa$-fold intersection property}.

\subsection{Uniformities and
homomorphisms}

 If $\mathcal U$ is a relation on an algebra $A$,
and $f:B
\to A$ is a homomorphism from another algebra $B$, then
we denote by $f^{-1}(U)$ the relation
$\{\,\pair b{b'}\mid f(b)\mathrel Uf(b')\,\}$. If
$\mathcal U$ is a filter of relations on $A$, then we
denote by $f^{-1}(\mathcal U)$ the filter
$\Fg\{\,f^{-1}(U)\mid U\in\mathcal U\,\}$.

\begin{proposition} Let $A$, $B$ be algebras, and $f:B\to
A$ a homomorphism. We have
\begin{enumerate}
\item[(1)] If $\mathcal U$ is a compatible uniformity
(compatible semiuniformity, singly compatible
semiuniformity) on
$A$, then
$f^{-1}(\mathcal U)$ is a compatible uniformity
(respectively, compatible semiuniformity, singly
compatible semiuniformity) on
$B$.
\item[(2)] The mapping $\mathcal U\mapsto f^{-1}(\mathcal
U)$ preserves arbitrary meets.
\end{enumerate}
\end{proposition}

Now, suppose that we have
$\mathcal U\in\Fil\Rel A$,
and a homomorphism
$f:A\to B$. If $\mathcal U\in\OpUnif A$, we define
$f_{*c}(\mathcal U)$ to be the meet of all $\mathcal
V\in\OpUnif B$ such that $\mathcal U\leq f^{-1}(\mathcal
V)$.

\begin{proposition} Let $A$ and $B$ be algebras, and
$f:A\to B$ a homomorphism. Then the mapping $\mathcal
U\mapsto f_{*c}(\mathcal U)$ preserves arbitrary joins.
\end{proposition}

Let $A$, $B$ be algebras, and $f:A\to B$ a homomorphism,
and $U$ a relation on $A$. If $n$, $n'\geq 0$ and $t$ is
an $(n+n')$-ary term, we denote by $L_{f,n,n',t}(U)$ the
set of pairs $\pair{t(\mathbf b,f(\mathbf a))}{t(\mathbf
b,f(\mathbf a'))}$ where $\mathbf b\in B^n$ and $\mathbf
a$,
$\mathbf a'\in A^{n'}$ are such that $\mathbf a\mathrel
U\mathbf a'$.

\begin{theorem}\label{T:Push} Given $A$, $B$, $f:A\to B$,
and a uniformity $\mathcal U$ on $A$, we have
\begin{enumerate}
\item[(1)] Given $n$, $n'$, and $t$, the set of relations
$L_{f,n,n',t}(U)$ on
$B$, for $U\in\mathcal U$, is a base for a filter
$\mathcal L_{f,n,n',t}(\mathcal U)$ of relations on $B$.
\item[(2)] Given $n$, $n'$, and $t$, if $\mathcal U$ is
compatible, then
$\mathcal L_{f,n,n',t}(\mathcal U)\leq f_{*c}(\mathcal
U)$.
\item[(3)] $\mathcal L_f(\mathcal
U)\overset{def}=\bigcap_{n,n',t}\mathcal L_{f,n,n',t}(U)$
is a singly compatible semiuniformity.
\item[(4)] $f_{*c}(\mathcal U)=\Ug(\mathcal L_f(\mathcal
U))$.
\end{enumerate}
\end{theorem}

\begin{proof} (1) is clear.

(2): Let $V\in f_{*c}(\mathcal U)$. There is a $\bar V\in
f_{*c}(\mathcal U)$ such that $t(\bar V)\subseteq V$.
There is a $U\in\mathcal U$ such that $U\subseteq
f^{-1}(\bar V)$, because $\mathcal U\leq
f^{-1}(f_{*c}(\mathcal U))$. Then
$f(U)\subseteq\bar V$ and $L_{f,n,n',t}(U)\subseteq V$.
Thus, $\mathcal L_{f,n,n',t}(\mathcal U)\leq
f_{*c}(\mathcal U)$.

(3): (U1) and (U2) are clear. If $t$ is the unary term
$t(x)=x$, then $\mathcal L_{f,1,0,t}(\mathcal
U)=\{\,\Delta_B\,\}$, proving (U3).
$L_{f,n,n',t}(U^{-1})=L_{f,n,n',t}(U)^{-1}$, so if
$\mathcal U$ is a semiuniformity, then so is $\mathcal
L_{f,n,n',t}(\mathcal U)$ for all $n$, $n'$, and $t$,
proving (U4). To show single compatibility, we must
show that given
$V\in\mathcal L_f(\mathcal U)$ and an $(\bar
n+1)$-ary term $\bar t$ for $\bar n\geq 0$, there is a
$\bar V\in\mathcal L_f(\mathcal U)$ such  that $\bar
t(\bar V,\mathbf c)\subseteq V$ for any $\mathbf c\in B
^{\bar n}$. It suffices to show that, given $n$, $n'$,
and $t$, there is a $U_{n,n',t}\in\mathcal U$ such that
$\bar t(L_{f,n,n',t}(U_{n,n',t}),\mathbf c)\subseteq V$
for $\mathbf c\in B^{\bar n}$. For, then $\hat
U=\bigcup_{n,n',t}L_{f,n,n',t}(U_{n,n',t})\in\mathcal
L_f(\mathcal U)$ and
$\bar t(\hat U,\mathbf c)\subseteq V$ for all $\mathbf c$.
Now, $\bar t(L_{f,n,n',t}(U),\mathbf c)\subseteq
L_{f,n+\bar n,n',\tilde t}(U)$ where $\tilde t(\mathrm
x,\mathbf y,\mathbf z)=\bar t(t(\mathrm x,\mathbf
z),\mathbf y)$. Thus, if we pick $U_{n,n',t}$ such that
$L_{f,n+\bar n,n',\tilde t}(U_{n,n',t})\subseteq V$, which
we can do because $\mathcal L_{f,n+\bar n,n',\tilde
t}(\mathcal U)\leq f_{*c}(\mathcal U)$ by (2), then we
have $\bar t(L_{f,n,n',t}(U_{n,n',t}),\mathbf c)\subseteq
V$.

(4): We have $\mathcal L_f(\mathcal U)\leq f_{*c}(\mathcal
U)$ by (2). Thus, $\Ug\mathcal L_f(\mathcal U)\leq
f_{*c}(\mathcal U)$. To show the opposite inequality, it
suffices to note that $\mathcal U\leq
f^{-1}(\mathcal L_f(\mathcal U))$. For, $\mathcal U\leq
f^{-1}(\mathcal L_{1,0,x}(\mathcal U))$.
\end{proof}

\begin{remark} In \cite[Section 11]{R02}, there is another
incorrect discussion about the procedure for finding the
colimit of a diagram $F:\mathcal D\to\mathbf V(\Unif)$ in
the category
$\mathbf V(\Unif)$. The uniformity of the colimit
is the smallest
\emph{compatible} uniformity greater than or equal to all
of the $\iota_d(\mathcal U(d))$, where $\iota_d$ is the
insertion of $F(d)$ into the colimit and $\mathcal U(d)$
is the uniformity on $F(d)$.
\end{remark}

\subsection{Completion} The \emph{completion} of an
algebra $A$ with respect to a compatible uniformity
$\mathcal U$ is defined as the set of equivalence classes
of Cauchy nets in $A$ with respect to $\mathcal U$, and
we denote it by $A/\mathcal U$. $A/\mathcal U$ has a
natural structure of uniform universal algebra. The
mapping from $A$ to $A/\mathcal U$ taking $a\in A$ to the
equivalence class containing the constant nets at $a$ is
a uniform homomorphism onto a dense subset of $A/\mathcal
U$, and we denote this mapping by $\nat\mathcal U$. Note
that in \cite{R02}, we used the notation $\eta_{\mathcal
U}$ for $\nat\mathcal U$.

We denote the natural uniformity on $A/\mathcal U$ by
$\mathcal U/\mathcal U$; it has a base of relations
$R(\mathcal U,U)$ for $U\in\mathcal U$, where $R(\mathcal
U,U)$ relates two equivalence classes $k$, $k'$ of Cauchy
nets iff there exist nets $n\in k$, $n'\in k'$ such that
$n(d)\mathrel Un'(d')$ for large enough $d$ and $d'$.

More generally, if $\mathcal V$ is another compatible
uniformity on $A$ such that $\mathcal U\leq\mathcal V$,
then there is a uniformity $\mathcal V/\mathcal U$ on
$A/\mathcal U$ having a base of relations $R(\mathcal
U,V)$ for $V\in\mathcal V$.

\begin{remark} It follows from \cite[Theorem
9.9(2)]{R02} that $\mathcal V/\mathcal U=(\nat\mathcal
U)_{*c}(\mathcal V)=(\nat\mathcal U)_*(\mathcal V)$.
\end{remark}

Formation of the completion plays the same role in
uniform universal algebra that formation of quotient
algebras plays in standard universal algebra, and
satisfies many of the same properties. See
\cite[Section~9]{R02}.
As an example of the close relationship between these two
constructions,
if $\alpha\in\Con A$ then $A/\Ug\{\,\alpha\,\}\cong
\pair{A/\alpha}{\mathcal U_d}$, where $\mathcal
U_d=\Fg\{\,\Delta\,\}$ is the discrete uniformity.

\subsection{Joins of compatible uniformities}

Joins in the lattice of compatible uniformities
are the same as in the lattice of uniformities on the
underlying set \cite[Theorem 5.3]{R02}.

\subsection{Compatible uniformities on algebras in
congruence-permutable algebras}

We recall (\cite[Theorems~6.4 and 6.2]{R02}) that if $A$
is an algebra
in a congruence-permutable variety, and $\mathcal U$,
$\mathcal V\in\OpUnif A$, then $\mathcal U$ and $\mathcal
V$ permute, and
that
$\OpUnif A$ is modular.

\section{Topological Groups and Uniform
Groups}\label{S:Groups}

A \emph{topological group} is a group object in the
category of topological spaces and continuous functions.
Such an object is determined by a group $G$ and a
topology $\mathcal T$ on $G$ such that the group
operations are continuous functions. Thus, it is
different from a uniform group, or group with a
compatible uniformity $\pair G{\mathcal U}$, where the
operations are required to be not only
continuous but uniformly continuous.

\subsection{Axioms for topological groups}

If $\pair G{\mathcal T}$ is a topological group, then the
set $\mathcal N$ of neighborhoods of the identity $e$
satisfies the following axioms:

\begin{enumerate}
\item[(G1)]
If $N\in\mathcal N$ and $N\subseteq N'$, then
$N'\in\mathcal N$.
\item[(G2)]
If $N$, $N'\in\mathcal N$, then $N\cap N'\in\mathcal N$.
\item[(G3)]
For any $N\in\mathcal N$, there exists a neighborhood
$\bar N\in\mathcal N$ such that $\bar N\bar N=\{\,xy\mid
x,y\in\bar N\,\}\subseteq N$.
\item[(G4)]
If $N\in\mathcal N$, then $N^{-1}\in\mathcal N$.
\item[(G5)]
If $N\in\mathcal N$ and $a\in G$, then
$a^{-1}Na\in\mathcal N$.
\end{enumerate}

A stronger version of (G5) which will be useful to us
is
\begin{enumerate}
\item[(G$5'$)] If $N\in\mathcal N$, then
$\bigcap_aa^{-1}Na\in\mathcal N$.
\end{enumerate}

\subsection{Uniform groups}

A \emph{uniform group}, or, \emph{group with a compatible
uniformity}, is a pair $\pair G{\mathcal U}$ where $G$ is
a group and $\mathcal U$ is a compatible uniformity. If
$\pair G{\mathcal U}$ is a uniform group, then the
topology $\mathcal T$ underlying the uniformity $\mathcal
U$ is compatible with the group operations. The
neighborhood system $\mathcal N$ of this topology is
given by $N\in\mathcal N$ iff $N=\{\,x\mid x\mathrel
Ue\,\}$ for some $U\in\mathcal U$. The conditions (G1)
through (G5) can easily be verified.

\subsection{Translation invariance}

If $A$ is a group, we say
that a relation
$U\subseteq A^2$ is \emph{left translation invariant}
(\emph{right translation invariant}) if
 $a\in A$ and $b\mathrel Ub'$ imply $ab\mathrel
Uab'$ (respectively, $ba\mathrel Ub'a$).
If $A$ is abelian, then left translation invariance and
right translation invariance coincide and we simply say
that a relation is \emph{translation invariant}.

\begin{lemma}\label{T:InvBase} Every compatible uniformity
$\mathcal U$ on a group $A$ has a base of
left translation invariant relations and a base of right
translation invariant relations.
\end{lemma}

\begin{proof}
Given $U\in\mathcal U$, let $U'\in\mathcal U$ be such
that if $b\mathrel{U'}c$, then $ab\mathrel U ac$ for
any $a$. Define the relation $V$ by $x\mathrel Vy$ iff
there exist $b$, $c$ such that $x^{-1}y=b^{-1}c$ and
$b\mathrel{U'}c$. Then clearly $U'\subseteq V$, so that
$V\in\mathcal U$, and also
$V\subseteq U$ because $b\mathrel{U'}c$ and
$x^{-1}y=b^{-1}c$ imply $x=(xb^{-1})b\mathrel
U(xb^{-1})c=xx^{-1}y=y$. But
$V$ is left translation invariant. This proves that
$\mathcal U$ has a base of left translation invariant
relations; the proof that $\mathcal U$ has a base of
right translation invariant relations is similar.
\end{proof}

\subsection{The left uniformity and right uniformity of a
compatible topology}

Let $\mathcal N$ be a neighborhood system for a
compatible topology on $G$. If $N\in\mathcal N$, we
define $N_l=\{\,\pair xy\in G^2\mid y\in xN\,\}$. The set
of relations $\{\,N_l\mid N\in\mathcal N\,\}$ is a base
for a compatible uniformity $\mathcal U_{\mathcal T,l}$ on
$G$, called the
\emph{left uniformity}.  Similarly, if $N\in\mathcal N$,
we define $N_r=\{\,\pair xy\mid y\in Nx\,\}$, and the
$N_r$ form a base for the \emph{right uniformity},
denoted by $\mathcal U_{\mathcal T,r}$.
 The inverse operation is
a uniform isomorphism (of the uniform structure) when
viewed as a function from
$\pair G{\mathcal U_{\mathcal T,l}}$ to
$\pair G{\mathcal U_{\mathcal T,r}}$.
Note that $\mathcal U_{\mathcal T,l}$ has a base of left
translation invariant relations, and $\mathcal
U_{\mathcal T,r}$ has a base of right translation
invariant relations.

\begin{theorem}\label{T:OnlyOne} Let $\pair G{\mathcal T}$
be a topological group. There is at most one compatible
uniformity $\mathcal U$ on $G$ such that $\mathcal T$ is
the topology underlying $\mathcal U$, and in this case,
$\mathcal U=\mathcal U_{\mathcal T,l}=\mathcal
U_{\mathcal T,r}$.
\end{theorem}

\begin{proof} If $\mathcal U$ exists, then by
Lemma~\ref{T:InvBase},
$\mathcal U$ has a base of left-invariant relations.
It follows that if $\mathcal T$ is the underlying topology
of a compatible uniformity $\mathcal U$, then $\mathcal
T$ determines $\mathcal U$ as $\mathcal U_{\mathcal T,l}$.
Similar arguments apply to the right uniformity.
\end{proof}

The Theorem causes us to ask the question: If we have a
topological group $\pair A{\mathcal T}$ such  that
$\mathcal U_{\mathcal T,l}=\mathcal U_{\mathcal T,r}$,
must this uniformity be compatible?

\begin{theorem} Let $\pair A{\mathcal T}$ be a
topological group. The following are equivalent:
\begin{enumerate}
\item[(1)] $\mathcal U_{\mathcal T,r}\leq\mathcal
U_{\mathcal T,l}$
\item[(2)] $\mathcal U_{\mathcal T,l}\leq\mathcal
U_{\mathcal T,r}$
\item[(3)] $\mathcal U_{\mathcal T,l}=\mathcal
U_{\mathcal T,r}$
\item[(4)] $\mathcal U_{\mathcal T,l}$ is compatible
\item[(5)] $\mathcal U_{\mathcal T,r}$ is compatible
\item[(6)] The neighborhood system $\mathcal N$
corresponding to $\mathcal T$ satisfies (G6${\,}'$).
\end{enumerate}
\end{theorem}

\begin{proof}
Clearly, (3)$\implies$(1) and (3)$\implies$(2).
We have (4)$\implies$(3) and (5)$\implies$(3) by
Theorem~\ref{T:OnlyOne}.

(1)$\implies$(6): If $\mathcal U_{\mathcal
T,r}\leq\mathcal U_{\mathcal T,l}$, then given
$N\in\mathcal N$, there is an
$\bar N\in\mathcal N$ such that $\bar N_r\subseteq N_l$.
That is,
$y\in\bar Nx\implies y\in xN$, or $yx^{-1}\in\bar
N\implies x^{-1}y\in N$. This implies by a change of
variables that
$x\in\bar N\implies a^{-1}xa\in N$ for all $x$ and $a$,
i.e., (6).

To prove (6)$\implies$(4), we must show that $\mathcal
U_{\mathcal T,l}$ is compatible with respect to group
multiplication and the inverse operation.

To show $\mathcal U_{\mathcal T,l}$ is compatible with
respect to group multiplication, it suffices to show that
given $N\in\mathcal N$, there is an $\bar N\in\mathcal N$
such that $\pair x{x'}$, $\pair y{y'}\in\bar
N_l\implies\pair{xy}{x'y'}\in N_l$, or in other words,
$x^{-1}x'$, $y^{-1}y'\in\bar N\implies
y^{-1}x^{-1}x'y'\in N$. Given $N$, there is an $\hat
N\in\mathcal N$ such that $a^{-1}\hat Na\subseteq N$ for
all $a$. There is an $\tilde N\in \mathcal N$ such that
$\tilde N\tilde N\subseteq\hat N$, by (G5'). Finally,
there is an
$\bar N\in\mathcal N$ such that $a^{-1}\bar
Na\subseteq\tilde N$ for all $a$. Then
\begin{align*}
x^{-1}x',\,y^{-1}y'\in\bar N&\implies
x^{-1}x',\,y'y^{-1}\in\tilde N\\
&\implies x^{-1}x'y'y^{-1}\in\hat N\\
&\implies y^{-1}x^{-1}x'y'\in N.
\end{align*}

To show that $\mathcal U_{\mathcal T,l}$ is compatible
with respect to the inverse operation, it suffices to
show that for each $N$, there is an $\bar N$ such  that
$\pair xy\in\bar N_l\implies\pair{x^{-1}}{y^{-1}}\in
N_l$, or in other words, $x^{-1}y\in\bar N\implies
xy^{-1}\in N$. Given $N$, let $\hat N$ be such that $\hat
N^{-1}\subseteq N$, and let $\bar N$ be such that
$xy\in\bar N\implies yx\in\bar N$ (true by (G$5'$)). Then
$x^{-1}y\in\bar N\implies yx^{-1}\in\hat N\implies
xy^{-1}\in N$.
Thus, (6)$\implies$(4).

The proof that (6)$\implies$(5)
is similar.
\end{proof}

\section{Term Conditions, Centralization, and Related
Commutator Operations}\label{S:Term}

In this section, we will discuss various conditions we
call \emph{term conditions}, and define notions of
centralization and commutator operations based on them.
First, more or less following \cite{F-McK}, we review the
\emph{$\alpha,\beta$-term condition}.
Then, we generalize this to uniformities and give the
\emph{$\mathcal U,\mathcal V$-term condition}. We also
give two weaker conditions, which we call the \emph{weak
$\alpha,\beta$-term condition} and the \emph{weak
$\mathcal U,\mathcal V$-term condition}. As we state
these four conditions, we give corresponding notions
of centralization. Then we define notions of commutator
(binary operations on $\Con A$ and $\OpUnif A$) derived
from the four types of centralization, and finally, we
prove some relationships between the various notions,
showing that centralization for congruences can be
considered a special case of centralization for
uniformities.

\subsection{The $\alpha,\beta$-term condition}

We begin with the
$\alpha,\beta$-term condition. We consider it as coming in
two equivalent forms:

If $\alpha$, $\beta$,
$\delta\in\Con A$, for some algebra $A$, then we say that
$\delta$ satisfies the \emph{first form of the
$\alpha,\beta$-term condition} if for all $a$, $a'\in A$
such that $a\mathrel\alpha a'$, for all $\mathbf b$,
$\mathbf c\in A^n$, $n>0$, such that $\mathbf
b\mathrel\beta\mathbf c$ (i.e., $b_i\mathrel\beta c_i$
for all $i$), and all $(n+1)$-ary terms $t$, $t(a,\mathbf
b)\mathrel\delta t(a,\mathbf c)$ implies $t(a',\mathbf
b)\mathrel\delta t(a',\mathbf c)$.

To give the second form of the $\alpha,\beta$-term
condition, we first define, given $n\geq 0$, $n'\geq 0$,
and
$t$, an $(n+n')$-ary term, and binary relations $U$, $V$
on
$A$, the set of $2\times 2$ matrices
\[
\mathrm M_{n,n,t}(U,V)=\left\{\left(\begin{array}{cc}
t(\mathbf a,\mathbf b) & t(\mathbf a,\mathbf c) \\
t(\mathbf a',\mathbf b) & t(\mathbf a',\mathbf c)
\end{array}\right)
\mathrel{\bigg|} \mathbf a, \mathbf a'\in
A^n,\,\mathbf
b,\mathbf c\in A^{n'},\,\text{ with }\mathbf a\mathrel
U\mathbf a'\text{ and }\mathbf b\mathrel V\mathbf
c\,\right\}.
\]
Then we say that $\delta$ satisfies the \emph{second form
of the $\alpha,\beta$-term condition} if for all $n$,
$n'$, and $t$,
$\left(\begin{array}{cc} u_{11} & u_{12} \\
u_{21} & u_{22}
\end{array}\right)\in\mathrm M_{n,n',t}(\alpha,\beta)$ and
$u_{11}\mathrel\delta u_{12}$ imply $u_{21}\mathrel\delta
u_{22}$.

\begin{proposition}\label{T:abTerm} The two forms of the
$\alpha,\beta$-term condition are equivalent.
\end{proposition}

\begin{proof} Clearly if $\delta$ satisfies the second
form, it satisfies the first. Given $\delta$ satisfying
the first form,
$n$,
$n'$, and $t$, and $\left(\begin{array}{cc} t(\mathbf
a,\mathbf b) & t(\mathbf a,\mathbf c) \\
t(\mathbf a',\mathbf b) & t(\mathbf a',\mathbf c)
\end{array}\right)\in\mathrm M_{n,n',t}(\alpha,\beta)$,
such that
$t(\mathbf a,\mathbf b)\mathrel\delta t(\mathbf a,\mathbf
c)$, we apply the first form
$n$ times, changing one component of $\mathbf a$ at a
time, to obtain the conclusion that $t(\mathbf a',\mathbf
b)\mathrel\delta t(\mathbf a',\mathbf c)$.
\end{proof}

In view of this equivalence, we simply say that $\delta$
\emph{satisfies the $\alpha,\beta$-term condition}.
When this is so, we also say that $\alpha$
\emph{centralizes} $\beta$ \emph{modulo} $\delta$, or
$C(\alpha,\beta;\delta)$.

\begin{proposition}\label{T:Cent} We have
\begin{enumerate}
\item[(1)] $C(\alpha,\beta;\alpha)$
\item[(2)] $C(\alpha,\beta;\beta)$
\item[(3)] $C(\alpha,\beta;\delta)$ and
$\alpha'\leq\alpha$,
$\beta'\leq\beta$ imply $C(\alpha',\beta';\delta)$
\item[(4)] If $C(\alpha,\beta;\delta_i)$ for $i\in I$,
then
$C(\alpha,\beta;\bigwedge_i\delta_i)$.
\end{enumerate}
\end{proposition}

\subsection{The $\mathcal U,\mathcal V$-term condition}
We now generalize the term condition to compatible
uniformities. Let $A$ be an algebra, and let $\mathcal U$,
$\mathcal V$, $\mathcal W\in\OpUnif A$.

We say that
$\mathcal W$ satisfies the \emph{first form of the
$\mathcal U,\mathcal V$-term condition} if for all $n>0$,
all
$(n+1)$-ary terms $t$, and all
$W\in\mathcal W$,
there are $U\in\mathcal U$, $V\in\mathcal V$, $\bar
W\in\mathcal W$ such that for all
$a$, $a'\in A$ such that $a\mathrel Ua'$, and all $\mathbf
b$,
$\mathbf c\in A^n$ such that $\mathbf b\mathrel
V\mathbf c$ for all
$i$,
$t(a,\mathbf b)\mathrel{\bar W}t(a,\mathbf c)$ implies
$t(a',\mathbf b)\mathrel Wt(a',\mathbf c)$.

We say that $\mathcal W$
satisfies the
\emph{second form of the $\mathcal U,\mathcal V$-term
condition} if for all
$n$, $n'$, all $(n+n')$-ary
$t$, and all $W\in\mathcal W$, there are $U\in\mathcal U$,
$V\in\mathcal V$, and $\bar W\in\mathcal W$ such that
$\left(\begin{array}{cc} u_{11} & u_{12} \\ u_{21} &
u_{22}
\end{array}\right)\in\mathrm M_{n,n',t}(U,V)$ and
$u_{11}\mathrel{\bar W}u_{12}$ imply $u_{21}\mathrel
Wu_{22}$.

\begin{theorem}\label{T:UVTerm} The two forms of the
$\mathcal U,\mathcal V$-term condition are equivalent.
\end{theorem}

\begin{proof} Clearly,
the second form of the condition implies the
first form. In the other
direction, given $n$, $n'$, and $t$, $t(\mathbf a,\mathbf
b)$ and
$t(\mathbf a,\mathbf c)$ can be changed to
$t(\mathbf{a'},\mathbf b)$ and $t(\mathbf{a'},\mathbf c)$
one component of $\mathbf a$ at a time.  For the
 $\alpha,\beta$-term condition, the result is immediate,
but in the uniformity-theoretic case, given $W\in\mathcal
W$, we must choose in reverse order $W_i$, $i=1$,
$\ldots,n$ such that the changes are valid, starting by
choosing $U_n\in\mathcal U$,
$V_n\in\mathcal V$, and $W_n\in\mathcal W$ such that
$a_n\mathrel{U_n}$, $\mathbf b\mathrel{V_n}\mathbf c$, and
$t(\langle a'_1,\ldots,a'_{n-1},a_n\rangle,\mathbf
b)\mathrel{W_n}t(\langle
a'_1,\ldots,a'_{n-1},a_n\rangle,\mathbf c)$ imply
$t(\mathbf a,\mathbf b)\mathrel W t(\mathbf{a'},\mathbf
c)$, and ending by choosing $U_1$, $V_1$,
and $W_1$ such that $a_1\mathrel{U_1}a'_1$, $\mathbf
b\mathrel{V_1}\mathbf c$, and
$t(\mathbf a,\mathbf b)\mathrel{W_1}t(\mathbf a,\mathbf
c)$ imply
$t(\langle
a'_1,a_2,\ldots,a_n\rangle,\mathbf
b)\mathrel{W_2}
t(\langle
a'_1,a_2,\ldots,a_n\rangle,\mathbf c)$. We then
let
$U=\bigcap_iU_i$, $V=\bigcap_iV_i$, and $\bar W=W_1$.
\end{proof}

If $\mathcal W$ satisfies the two equivalent forms of the
$\mathcal U,\mathcal V$-term condition, then
we say that
$\mathcal U$
\emph{centralizes}
$\mathcal V$
\emph{modulo} $\mathcal W$, or $C(\mathcal U,\mathcal
V;\mathcal W)$.

\begin{proposition}\label{T:UCent} We have
\begin{enumerate}
\item[(1)] $C(\mathcal U,\mathcal V;\mathcal U)$
\item[(2)] $C(\mathcal U,\mathcal V;\mathcal V)$
\item[(3)] $C(\mathcal U,\mathcal V;\mathcal W)$ and
$\mathcal U'\leq\mathcal U$, $\mathcal V'\leq\mathcal V$
imply $C(\mathcal U',\mathcal V';\mathcal W)$
\item[(4)] If $C(\mathcal U,\mathcal V;\mathcal W_i)$ for
$i\in I$, then $C(\mathcal U,\mathcal
V;\bigwedge_i\mathcal W_i)$
\item[(5)]
$C(\Ug\{\,\alpha\,\},\Ug\{\,\beta\,\};\Ug\{\,\delta
\,\})$ iff $C(\alpha,\beta;\delta)$.
\end{enumerate}
\end{proposition}

\begin{proof} To
show (1), let
$U\in\mathcal U$ and $t$ be given.
There is a symmetric $\bar U\in\mathcal U$ such that
$a\mathrel{\bar U}a'$ implies $t(a,\mathbf
b)\mathrel{{^3U}}t(a',\mathbf b)$ and
$t(a,\mathbf
c)\mathrel{{^3U}}t(a',\mathbf c)$. Then for any
$V\in\mathcal V$,
$a\mathrel{\bar
U}a'$ and
$t(a,\mathbf
b)\mathrel{{^3U}}t(a,\mathbf
c)$ imply $t(a',\mathbf b)\mathrel
Ut(a',\mathbf c)$.

To show (2), let $V\in\mathcal V$ and
$t$ be given. Then for some $\bar
V\in\mathcal V$,
$\mathbf b\mathrel V\mathbf c$ implies $t(a',\mathbf
b)\mathrel V t(a',\mathbf c)$ for any $a'$, by uniform
continuity of $t$, regardless of any consideration of
$t(a,\mathbf b)$ and $t(a,\mathbf c)$.

(3) is obvious.

(4): Suppose the $\mathcal W_i$ satisfy the $\mathcal
U$,$\mathcal V$-term condition, and that $n>0$, an
$(n+1)$-ary term $t$, and $W\in\bigwedge_i\mathcal W_i$
are given. Then $\bigwedge_{j=1}^kW_j\subseteq W$ for
some uniform neighborhoods $W_j\in W_{i_j}$, $i_j$ being
selected values of the index $i$. It suffices to show that
$\bigwedge_{j=1}^k\mathcal W_{i_j}$ satisfies the
$\mathcal U$,$\mathcal V$-term condition.

Let $U_j\in\mathcal U$, $V_j\in\mathcal V$, $\bar
W_j\in\mathcal W_{i_j}$ be relations as promised by the
$\mathcal U$,$\mathcal V$-term condition for $\mathcal
W_{i_j}$, and let $U=\bigwedge_jU_j$, $V=\bigwedge_jV_j$,
and
$\bar W=\bigwedge_jW_j$. Then if $a \mathrel U a'$,
$\mathbf b\mathrel V\mathbf c$, and $t(a,\mathbf
b)\mathrel{\bar W}t(a,\mathbf c)$, we have $t(a',\mathbf
b)\mathrel{W_j}t(a',\mathbf c)$ for each $j$, whence
$t(a',\mathbf b)\mathrel Wt(a',\mathbf c)$.

(5) follows easily from the definitions.
\end{proof}

\subsection{The weak term conditions}
We say that $\delta$ satisfies
the \emph{weak $\alpha,\beta$-term condition} if for all
$n$, $n'$, and $t$, $\left(\begin{array}{cc}
u_{11} & u_{12} \\
u_{21} & u_{22}\end{array}\right)\in
\mathrm M_{n,n',t}(\alpha,\beta)$ and $u_{11}=u_{12}$
imply
$u_{21}\mathrel\delta u_{22}$.  In case $\delta$
satisfies the weak $\alpha,\beta$-term condition, we say
that $\alpha$ \emph{weakly centralizes}
$\beta$ \emph{modulo} $\delta$, or $\tilde
C(\alpha,\beta;\delta)$.

\begin{proposition}\label{T:CentW} We have
\begin{enumerate}
\item[(1)] $\tilde C(\alpha,\beta;\alpha)$
\item[(2)] $\tilde C(\alpha,\beta;\beta)$
\item[(3)] $\tilde C(\alpha,\beta;\delta)$ and
$\alpha'\leq\alpha$, $\beta'\leq\beta$ imply
$C(\alpha',\beta';\delta)$
\item[(4)] If $\tilde C(\alpha,\beta;\delta_i)$ for $i\in
I$, then $\tilde C(\alpha,\beta;\bigwedge_i\delta_i)$
\item[(5)] $\tilde C(\alpha,\beta;\delta)$ and
$\delta\leq\delta'$ imply $\tilde C(\alpha,\beta;\delta')$
\item[(6)] $C(\alpha,\beta;\delta)\implies\tilde
C(\alpha,\beta;\delta)$.
\end{enumerate}
\end{proposition}

We say that $\mathcal W$
satisfies the \emph{weak\/ $\mathcal U,\mathcal V$-term
condition} if for all $n$, $n'$, and $t$, and all
$W\in\mathcal W$, there exist $U\in\mathcal U$ and
$V\in\mathcal V$ such that
$\left(\begin{array}{cc}
u_{11} & u_{12} \\
u_{21} & u_{22}\end{array}\right)\in\mathrm M_{n,n',t}(
U,V)$ and $u_{11}=u_{12}$ imply $u_{21}\mathrel W u_{22}$.
If $\mathcal W$ satisfies the
weak $\mathcal U,\mathcal V$-term condition, we say that
$\mathcal W$ \emph{weakly centralizes} $\mathcal U$
\emph{modulo} $\mathcal V$, or $\tilde C(\mathcal
U,\mathcal V;\mathcal W)$.

\begin{proposition}\label{T:UCentW} We have
\begin{enumerate}
\item[(1)] $\tilde C(\mathcal U,\mathcal V;\mathcal U)$
\item[(2)] $\tilde C(\mathcal U,\mathcal V;\mathcal V)$
\item[(3)] $\tilde C(\mathcal U,\mathcal V;\mathcal W)$
and
$\mathcal U'\leq\mathcal U$, $\mathcal V'\leq\mathcal V$
imply $\tilde C(\mathcal U',\mathcal V';\mathcal W)$
\item[(4)] If $\tilde C(\mathcal U,\mathcal V;\mathcal
W_i)$ for $i\in I$, then $\tilde C(\mathcal U,\mathcal
V;\bigwedge_i\mathcal W_i)$
\item[(5)] $\tilde C(\mathcal U,\mathcal V;\mathcal W)$
and
$\mathcal W\leq\mathcal W'$ imply $\tilde C(\mathcal
U,\mathcal V;\mathcal W')$
\item[(6)] $C(\mathcal U,\mathcal V;\mathcal W)\implies
\tilde C(\mathcal U,\mathcal V;\mathcal W)$
\item[(7)] $\tilde C(\Ug\{\,\alpha\,\},\Ug\{\,\beta\,\};
\Ug\{\,\delta\,\})$ iff $\tilde C(\alpha,\beta;\delta)$.
\end{enumerate}
\end{proposition}

\begin{proof}
(1): Given $n$, $n'$, and $t$, and $U\in\mathcal U$, there
exists a symmetric $U'\in\mathcal U$ such that $\mathbf
a\mathrel{U'}\mathbf a'$ implies $t(\mathbf a,\mathbf
b)\mathrel{^2U}t(\mathbf a',\mathbf b)$ for all $\mathbf
b$. It follows that for any $V\in\mathcal V$,
$\left(\begin{array}{cc} u_{11}&u_{12}\\
u_{21}&u_{22}\end{array}\right)\in\mathrm
M_{n,n',t}(U',V)$ and
$u_{11}=u_{12}$ imply
$u_{21}\mathrel{^2U}u_{11}=u_{12}\mathrel{^2U}u_{22}$.

(2) The proof of $\tilde C(\mathcal U,\mathcal
V;\mathcal V)$ is the same as the proof of $C(\mathcal
U,\mathcal V;\mathcal V)$.

(3) is obvious.

(4): Again we can reduce to the case $I$ finite. For
$W\in\bigwedge_i\mathcal W_i$,
$W=\bigcap_iW_i$ for some $W_i\in\mathcal W_i$. Given
$n$, $n'$, and $t$, there exist $U_i\in\mathcal U$,
$V_i\in\mathcal V$ such that $\left(\begin{array}{cc}
u_{11}&u_{12}\\ u_{21}&u_{22}\end{array}\right)\in\mathrm
M_{n,n',t}(U_i,V_i)$ and $u_{11}=u_{12}$ imply
$u_{21}\mathrel{W_i}u_{22}$. Then
$\left(\begin{array}{cc} u_{11}&u_{12}\\
u_{21}&u_{22}\end{array}\right)\in
\mathrm M_{n,n',t}(\bigcap_iU_i,\bigcap_iV_i)$ and
$u_{11}=u_{12}$ imply $u_{21}\mathrel Wu_{22}$.

(5), (6), and (7) are clear.
\end{proof}

\subsection{Related commutator operations}

These notions of centralization lead to definitions for
binary operations on the lattice of compatible
uniformities.
Recall \cite{F-McK} that if $A$ is an algebra, and
$\alpha$,
$\beta\in\Con A$, then
$C(\alpha,\beta)$ is defined as the least congruence on
$A$ satisfying the
$\alpha,\beta$-term condition.
Similarly, we denote by
$C(\mathcal U,\mathcal V)$ the least uniformity
satisfying the $\mathcal U,\mathcal V$-term condition,
by
$\tilde C(\alpha,\beta)$ the least congruence
satisfying the weak $\alpha,\beta$-term condition, and by
$\tilde C(\mathcal U,\mathcal V)$ the least uniformity
satisfying the weak $\mathcal U,\mathcal V$-term
condition. 
These uniformities exist by Proposition~\ref{T:UCent}(4)
and Proposition~\ref{T:UCentW}(4).

These commutator operations have some common
properties:

\begin{theorem}\label{T:AllComm} Let $\hat C(x,y)$ stand
for
$C(\alpha,\beta)$, $C(\mathcal U,\mathcal V)$, $\tilde
C(\alpha,\beta)$, or $\tilde C(\mathcal U,\mathcal V)$.
Then
\begin{enumerate}
\item[(1)] $\hat C(x,y)\leq x\wedge y$
\item[(2)] $x'\leq x$, $y'\leq y$ imply $\hat
C(x',y')\leq\hat C(x,y)$.
\end{enumerate}
\end{theorem}

Here are some explicit formulas for $C(\alpha,\beta)$ and
$\tilde C(\alpha,\beta)$:

\begin{proposition}\label{T:Explic}
Let $A$ be an algebra, and $\alpha$, $\beta\in\Con A$.
We have
\begin{enumerate}
\item[(1)] $C(\alpha,\beta)=\bigcup_\nu R_\nu$, where the
relations
$R_\nu$ are defined for all ordinal numbers $\nu$, as
follows:
\[R_\nu=\begin{cases} \Delta^A, &\nu=0 \\
\Cg\{\,\pair {t(a',\mathbf b)}{t(a',\mathbf
c)}\mid t\text{ is $(n+1)$-ary, }a\mathrel\alpha
a',\,\mathbf b\mathrel\beta\mathbf
c,\\\qquad\qquad\qquad\text{ and }
 t(a,\mathbf
b)\mathrel{R_{\nu'}}t(a,\mathbf c)\,\},
& \nu=\nu'+1\\
\bigcup_{\nu'<\nu}R_{\nu'}, & \text{for $\nu$ a limit
ordinal}
\end{cases}
\]
\item[(2)] $\tilde
C(\alpha,\beta)=\Cg\tilde R$ where
\[
\tilde R=\{\,\pair{u_{21}}{u_{22}}\mathrel{\bigg|}
\left(\begin{array}{cc}u_{11}&u_{12}\\
u_{21}&u_{22}\end{array}\right)\in
\mathrm M_{n,n',t}(\alpha,\beta)\text{ for some $n$, $n'$,
and
$t$, and $u_{11}=u_{12}$}\,\}.
\]
\end{enumerate}
\end{proposition}

\begin{theorem}\label{T:CommCong}
If $\alpha$, $\beta\in\Con A$, then we have
\begin{enumerate}
\item[(1)] $
C(\Ug\{\,\alpha\,\},\Ug\{\,\beta\,\})=\Ug\{\,
C(\alpha,\beta)\,\}$
\item[(2)] $\tilde
C(\Ug\{\,\alpha\,\},\Ug\{\,\beta\,\})=\Ug\{\,\tilde
C(\alpha,
\beta)\,\}$.
\end{enumerate}
\end{theorem}

\begin{proof} (1): By Proposition~\ref{T:UCent}(5),
since
$C(\alpha,\beta)$ satisfies the $\alpha,\beta$-term
condition,
the uniformity $\Ug\{\,C(\alpha,\beta)\,\}$ satisfies the
$\Ug\{\,\alpha\,\},\Ug\{\,\beta\,\}$-term condition.
Thus,
$C(\Ug\{\,\alpha\,\},\Ug\{\,\beta\,\})\leq\Ug\{\,
C(\alpha,\beta)\,\}$.
To show the opposite inequality, we
must show that $W\in
C(\Ug\{\,\alpha\,\},\Ug\{\,\beta\,\})$ implies
$C(\alpha,\beta)\subseteq W$. By
Proposition~\ref{T:Explic}(1),
$C(\alpha,\beta)=\bigcup_\nu R_\nu$. However, by
transfinite induction, and the
$\Ug\{\,\alpha\,\},\Ug\{\,\beta\,\}$-term condition, we
also have $R_\nu\subseteq W$ for all $\nu$ and
$W$. Thus,
$C(\alpha,\beta)\subseteq
W$ and $\Ug\{\,C(\alpha,\beta)\,\}\leq
C(\Ug\{\,\alpha\,\},\Ug\{\,\beta\,\})$.

(2): The proof that $\tilde
C(\Ug\{\,\alpha\,\},\Ug\{\,\beta\,\})\leq
\Ug\{\,\tilde C(\alpha,\beta)\,\}$ follows that same
pattern as for the operation $C(-,-)$. To show that
$\Ug\{\,\tilde C(\alpha,\beta)\,\}\leq\tilde
C(\Ug\{\,\alpha\,\},\Ug\{\,\beta\,\})$, we must show that
for every $W\in\tilde
C(\Ug\{\,\alpha\,\},\Ug\{\,\beta\,\})$, $\tilde
C(\alpha,\beta)\subseteq W$.
However, it is clear that $\tilde R\in W$, where
$\tilde R$ is the relation defined in the statement of
Proposition~\ref{T:Explic}(2).
Thus,
$\tilde R\subseteq\bigcap\tilde
C(\Ug\{\,\alpha\,\},\Ug\{\,\beta\,\})$. But this
intersection is a congruence. Thus, $\tilde
C(\alpha,\beta)=\Cg\tilde R\subseteq W$ by
Proposition~\ref{T:Explic}(2).
\end{proof}

By this theorem, the commutator operations $C(-,-)$ and
$\tilde C(-,-)$ on uniformities extend the
corresponding commutator operations
on congruences, and we can compute the
commutators on congruences by computing the
corresponding commutators of uniformities. The rule is to
promote both arguments to uniformities, and then take the
chosen commutator. The resulting uniformity then
determines the desired congruence.

\section{The Commutator on Uniformities in
Congruence-Modular Varieties}\label{S:Commutator}

As described in the previous section, the situation for
a general variety is that we have defined two possibly
different, possibly noncommutative commutator operations
on uniformities, $C(-,-)$ and $\tilde C(-,-)$. We will
show in this section that as it is with congruences
\cite{F-McK}, the situation is much simplified for
congruence-modular varieties: these operations coincide
and are commutative.

\subsection{$\mathcal M(\mathcal U,\mathcal V)$,
$\mathrm x_{\mathbf m}(\mathcal M)$, and $\mathcal
X_{\mathbf m}(\mathcal U,\mathcal V)$}

\begin{proposition} Let $\mathcal U$, $\mathcal
V\in\OpUnif A$. The set of sets of
$2\times 2$ matrices
$\mathrm M_{n,n',t}(U,V)$, $U\in\mathcal
U$, $V\in\mathcal V$ is a base for a filter $\mathcal
M_{n,n',t}(\mathcal U,\mathcal V)$ of sets of $2\times
2$ matrices of elements of $A$.
\end{proposition}

If $\mathcal U$, $\mathcal V\in\OpUnif A$, then we
define
\[\mathcal M(\mathcal U,\mathcal
V)=\bigvee_{n,n',t}\mathcal M_{n,n',t}(\mathcal
U,\mathcal V)=\bigcap_{n,n',t}\mathcal
M_{n,n',t}(\mathcal U,\mathcal V).
\]

Let
$m_0$,
$m_1$,
$\ldots$,
$m_k$ be a sequence of quaternary terms for
$A$. If $M$ is a set of $2\times 2$-matrices of elements
of
$A$,
we denote by $\mathrm x_{\mathbf m}(M)$ the set of pairs
\[\pair{m_i(a,b,d,c)}{m_i(a,a,c,c)}
\]
such that $i\leq k$ and $\left(\begin{array}{cc}a&b\\
c&d\end{array}\right)\in M$.

\begin{proposition} Given a filter $\mathcal M$ of sets
of $2\times 2$-matrices of elements of $A$, the set of
$\mathrm x_{\mathbf m}(M)$, $M\in\mathcal M$ forms a base
for a filter $\mathrm x_{\mathbf m}(\mathcal M)$ of
relations on $A$.
\end{proposition}

If $n$, $n'$, and $t$ are given, then we denote by
$\mathcal X_{\mathbf m,n,n',t}(\mathcal U,\mathcal V)$
the filter $\mathrm x_{\mathbf m}(\mathcal
M_{n,n',t}(\mathcal U,\mathcal V))$, and by
$\mathcal X_{\mathbf m}(\mathcal U,\mathcal V)$ the
filter $\mathrm x_{\mathbf m}(\mathcal M(\mathcal
U,\mathcal V))$.

\subsection{The commutator in modular varieties}

The following theorem can be compared to
\cite[Proposition 4.2(1)]{F-McK}, which establishes the
analogous result for commutators of congruences.

\begin{theorem}\label{T:Commutator} If
$\mathcal U$, $\mathcal V$, $\mathcal W$ are compatible
uniformities on $A$, an algebra in a
congruence-modular variety with sequence of Day terms
$\mathbf m$, then the following are equivalent:
\begin{enumerate}
\item[(1)] $\mathcal X_{\mathbf
m}(\mathcal U,\mathcal V)\leq\mathcal W$
\item[(2)] $\mathcal X_{\mathbf
m}(\mathcal V,\mathcal U)\leq\mathcal W$
\item[(3)] $C(\mathcal U,\mathcal V;\mathcal W)$
\item[(4)] $C(\mathcal V,\mathcal U;\mathcal W)$
\item[(5)] $\tilde C(\mathcal U,\mathcal V;\mathcal W)$
\item[(6)] $\tilde C(\mathcal V,\mathcal U;\mathcal W)$.
\end{enumerate}
\end{theorem}

\begin{proof} Obviously $(3)\implies(5)$ and
$(4)\implies(6)$. We will show that also
$(5)\implies (1)\implies (4)$. Exchanging $\mathcal U$ and
$\mathcal V$ will then give $(6)\implies (2)\implies (3)$,
completing the proof of equivalence.

(5)$\implies$(1): Assume $\tilde C(\mathcal
U,\mathcal V;\mathcal W)$. It suffices to show that if
$n$, $n'$,
$t$, and
$W\in\mathcal W$ are given, then there are
$U_{n,n',t}\in\mathcal U$ and
$V_{n,n',t}\in\mathcal V$ such that for all $\mathbf a$,
$\mathbf{a'}\in A^n$ with
$\mathbf a\mathrel{U_{n,n',t}}\mathbf{a'}$, and all
$\mathbf b$,
$\mathbf c\in A^{n'}$ such that $\mathbf b\mathrel
{V_{n,n',t}}\mathbf c$, and for all $i$, we have
$x\mathrel W
y$, where
\[
x=x(\mathbf a,\mathbf{a'},\mathbf b,\mathbf c)
=m_i(t(\mathbf a,\mathbf b),t(\mathbf a,\mathbf c),
t(\mathbf{a'},\mathbf c),t(\mathbf{a'},\mathbf b)),
\]
and
\[
y=y(\mathbf a,\mathbf{a'},\mathbf
b,\mathbf c)
=m_i(t(\mathbf a,\mathbf b),t(\mathbf a,\mathbf b),
t(\mathbf{a'},\mathbf b),t(\mathbf{a'},\mathbf b)).
\]
For, this implies that $\mathrm x_{\mathbf m}(\mathrm
M_{n,n',t}(U_{n,n',t},V_{n,n',t})\subseteq W$, and we
then have
\[
\mathrm x_{\mathbf m}(\bigcup_{n,n',t}\mathrm
M_{n,n',t}(U_{n,n',t},V_{n,n',t}))\subseteq W;
\]
but
$\bigcup_{n,n',t}\mathrm
M_{n,n',t}(U_{n,n',t},V_{n,n',t})\in\mathcal M(\mathcal
U,\mathcal V)$. Thus, it follows that $\mathcal
M(\mathcal U,\mathcal V)\leq\mathcal W$.

If we replace $\mathbf a$ by
$\mathbf{a'}$ at its second occurrence in the right-hand
expressions for $x$ and $y$, and $\mathbf{a'}$ by $\mathbf
a$ in its second occurrence, then we obtain expressions
for the same element
$z=m_i(t(\mathbf a,\mathbf b),t(\mathbf a',\mathbf c),
t(\mathbf a',\mathbf c),t(\mathbf a,\mathbf b))
=m_i(t(\mathbf a,\mathbf b),t(\mathbf a',\mathbf b),
t(\mathbf a',\mathbf b),t(\mathbf a,\mathbf b))=
t(\mathbf a,\mathbf b)$.
Let $s_i$ be the $(4n+4n')$-ary term given by 
\[s(\mathbf
g^1,\mathbf g^2,\mathbf g^3,\mathbf g^4,\mathbf
h^1,\mathbf h^2,\mathbf h^3,\mathbf h^4) =m_i(t(\mathbf
g^1,\mathbf h^1),t(\mathbf g^2,\mathbf h^2), t(\mathbf
g^3,\mathbf h^3),t(\mathbf g^4,\mathbf h^4)).
\]
Since $\tilde C(\mathcal U,\mathcal V;\mathcal W)$, there
are
$U_i\in\mathcal U$ and $V_i\in\mathcal V$ such that for
all 
$\left(
\begin{array}{cc}u_{11}&u_{12}\\u_{21}&u_{22}\end{array}
\right)\in\mathrm M_{4n,4n',s}(U_i,V_i)$,
$u_{11}=u_{12}$ implies $u_{21}\mathrel Wu_{22}$. But
\[
\left(\begin{array}{cc}
z&
z\\
x&
y
\end{array}\right)
=\left(\begin{array}{cc}
s_i(\mathbf a,\mathbf{a'},\mathbf{a'},\mathbf a,\mathbf
b,\mathbf c,\mathbf c,\mathbf b)&
s_i(\mathbf a,\mathbf{a'},\mathbf{a'},\mathbf a,\mathbf
b,\mathbf b,\mathbf b,\mathbf b)\\
s_i(\mathbf a,\mathbf a,\mathbf{a'},\mathbf{a'},\mathbf
b,\mathbf c,\mathbf c,\mathbf b)&
s_i(\mathbf a,\mathbf a,\mathbf{a'},\mathbf{a'},\mathbf
b,\mathbf b,\mathbf b,\mathbf b)
\end{array}\right)
\in\mathrm M_{4n,4n',s_i}(U_i,V_i).
\]
Thus, $x\mathrel Wy$. Letting $U=\bigwedge_iU_i$ and
$V=\bigwedge_iV_i$, we have $\mathrm
x_{\mathbf m}(M_{n,n',t}(U,V))\subseteq W$, proving that
(5)$\implies$(1).

To prove (1)$\implies$(4), let $W\in\mathcal W$ and let
$n$,
$n'$, and $t$ be given.  By
\cite[Lemma 7.1]{R02} (a uniformity-theoretic
generalization of
\cite[Lemma 2.3]{F-McK}) there exists
$\bar W\in\mathcal W$ such that $a$, $b$, $c$, $d\in A$
with $b\mathrel{\bar W}d$ and $m_i(a,a,c,c)\mathrel{\bar
W}m_i(a,b,d,c)$ for all $i$ imply $a\mathrel Wc$.
By (1), there are $U$ and $V$
such that $\mathrm x_{\mathbf
m}(M_{n',n,t}(U^{-1},V^{-1}))\subseteq\bar W^{-1}$. If
$\left(\begin{array}{cc}u_{11}&u_{12}\\u_{21}&u_{22}
\end{array}\right)\in\mathrm M_{n,n',t}(V,U)$,
then
$\left(\begin{array}{cc}u_{21}&u_{11}\\u_{22}&u_{12}
\end{array}\right)\in\mathrm M_{n',n,t}(U^{-1},V^{-1})$
and
\[m_i(u_{21},u_{21},u_{22},u_{22})\mathrel{\bar
W}m_i(u_{21},u_{11},u_{12},u_{22})\]
for all $i$. It
follows that if $u_{11}\mathrel{\bar W}u_{12}$, then
$u_{21}\mathrel Wu_{22}$. A $\bar W$, $U$,
and $V$ exists for each
$n$,
$n'$, and $t$, implying (4).
\end{proof}

If $\mathcal U$, $\mathcal V\in\OpUnif A$, we define
$[\mathcal U,\mathcal V]$ to be the least
$\mathcal W$ such that the six equivalent statements in
the theorem hold. Of course, we then have $[\mathcal
U,\mathcal V]=C(\mathcal U,\mathcal V)=C(\mathcal
V,\mathcal U)=\tilde C(\mathcal U,\mathcal V)=\tilde
C(\mathcal V,\mathcal U)=\Ug\mathcal X_{\mathbf
m}(\mathcal U,\mathcal V)$.

\begin{corollary} Let $A$ be an algebra in  a
congruence-modular variety with Day terms
$\mathbf m$, and let $\mathcal U$, $\mathcal V$, $\mathcal
W$, and $\mathcal W'\in\OpUnif A$. If
$C(\mathcal U,\mathcal V;\mathcal W)$ and $\mathcal
W<\mathcal W'$, then
$C(\mathcal U,\mathcal V;\mathcal W')$.
\end{corollary}

\section{$\OpUnif A(\alpha)$ and the Commutator
$[\Ug\{\,\alpha\,\},\mathcal U]$}\label{S:A(alpha)}

Recall that if $A$ is an algebra, and $\alpha\in\Con A$,
then
$A(\alpha)$ is the subalgebra of $A^2$ of pairs $\pair
ab$ such that $a\mathrel\alpha b$. We will denote by
$\pi$, $\pi':A(\alpha)\to A$ and $\Delta_\alpha:A\to
A(\alpha)$ the homomorphisms defined respectively by
$\pair ab\mapsto a$, $\pair ab\mapsto b$, and
$a\mapsto\pair aa$. (Note that in \cite{F-McK}, the
notation $\Delta_A$ is used for $\Delta_\alpha$, whereas
we use $\Delta_A$ to denote the diagonal set in $A^2$.)

If $\alpha$, $\beta\in\Con A$, where $A$ is an algebra,
then we can construct a congruence
$\Delta_{\alpha,\beta}\in\Con A(\alpha)$ by extending
$\beta$ along $\Delta_\alpha$.  That is,
$\Delta_{\alpha,\beta}=\Cg\{\,\pair{\pair aa}{\pair
bb}\mid a\mathrel\beta b\,\}$.

\begin{theorem} If $A$ is an algebra in a
congruence-modular variety, then
$(\pi')^{-1}[\alpha,\beta]=(\Delta_{\alpha,\beta}\wedge
\ker\pi)\vee\ker\pi'$.
\end{theorem}

\begin{proof} See
\cite[Theorem~4.9 and Exercise~4.4]{F-McK}.
\end{proof}

 The Theorem gives one way to compute
the commutator
$[\alpha,\beta]$. It may help to
illuminate the relationship between $[\alpha,\beta]$ and
$\Con A(\alpha)$ to note that the interval sublattice
$\I_{\Con A(\alpha)}[\bot,\ker\pi]$ transposes up to the
interval sublattice
$\I_{\Con
A(\alpha)}[\ker\pi',\top]$, which is of course isomorphic
to $\Con A$.

In this section, we will try to duplicate this
result with $\beta$ replaced by a compatible uniformity
$\mathcal U$, and $\Con A(\alpha)$ replaced by $\OpUnif
A(\alpha)$, the lattice of compatible uniformities of
$A(\alpha)$.

For this section, we will write $[\alpha,\mathcal U]$,
$\mathcal M(\alpha,\mathcal U)$, etc.\ for
$[\Ug\{\,\alpha\,\},\mathcal U]$, $\mathcal
M(\Ug\{\,\alpha\,\},\mathcal U)$, etc.

\subsection{$\mathcal
M(\alpha,\mathcal U)$ and
$\Delta_{\alpha,\mathcal U}$}

We define $\Delta_{\alpha,\mathcal U}$ as the analog of
$\Delta_{\alpha,\beta}$, that is, the compatible extension
$(\Delta_\alpha)_{*c}(\mathcal U)$ along $\Delta_\alpha$
of
$\mathcal U\in\OpUnif A$ to $A(\alpha)$.
We will construct
$\Delta_{\alpha,\mathcal U}$ using the filter
$\mathcal M(\alpha,\mathcal U)$ defined in
Section~\ref{S:Commutator}.

Now, $\mathcal M(\alpha,\mathcal U)$ is a
filter of subsets of $2\times 2$ matrices. We view each
matrix
$\left(\begin{array}{cc}
a&b\\c&d\end{array}\right)$ as a pair
$\pair{\left(\begin{array}{c}a\\c
\end{array}\right)}
{\left(\begin{array}{c}b\\d
\end{array}\right)}$ of elements of $A(\alpha)$ written
as column vectors. Thus,
$\mathcal M(\alpha,\mathcal U)$ can be viewed
as a filter of binary relations on $A(\alpha)$.

\begin{theorem} Viewed in this way, $\mathcal
M(\alpha,\mathcal U)$ is a singly compatible
semiuniformity on
$A(\alpha)$, and $\Delta_{\alpha,\mathcal U}=\Ug\mathcal
M(\alpha,\mathcal U)$.
\end{theorem}

\begin{proof}
We have $\mathcal M(\alpha,\mathcal U)=\mathcal
L_{\Delta_{\alpha}}(\mathcal U)$. Thus, by
Theorem~\ref{T:Push}, $\mathcal M(\alpha,\mathcal U)$ is a
singly compatible semiuniformity and
$\Delta_{\alpha,\mathcal
U}=(\Delta_\alpha)_{*c}(\mathcal U)=\Ug(\mathcal
L_{\Delta_{\alpha}}(\mathcal U))=\Ug\mathcal
M(\alpha,\mathcal U)$.
\end{proof}

\subsection{$[\alpha,\mathcal U]'$ and
$[\alpha,\mathcal U]$}

So far, it has not been proved that the lattice $\OpUnif
A(\alpha)$ of compatible uniformities on $A(\alpha)$ is
modular when $A$ belongs to a congruence-modular variety.
However, it is true in the congruence-modular case that if
$\mathcal U\in\OpUnif A(\alpha)$, $\mathcal
U\leq\ker\pi$, then
\[\mathcal U=(\mathcal U\vee\ker\pi')\wedge\ker\pi.\]
For, the fact that $\ker\pi$ is a congruence causes the
modular law to be true in this special case, by
\cite[Theorem~7.5]{R02}.

Thus, the interval
$\I_{\OpUnif A(\alpha)}[\bot,\ker\pi]$ can be embedded
into $\OpUnif A$ via a mapping $\phi$ such
that
$(\pi')^{-1}(\phi(\mathcal U))=\mathcal U\vee\ker\pi'$. We
shall see that if
$\mathcal U\in\OpUnif A$, where $A$ is an algebra with an
underlying group structure, then
$[\alpha,\mathcal U]$ does belong to the image of $\phi$,
and
$[\alpha,\mathcal U]=\phi(\Delta_{\alpha,\mathcal
U}\wedge\ker\pi)$. (Note that in this case, the modular
law also holds in $\OpUnif A$, by
\cite[Theorem~6.4]{R02}.)

Let us define
$[\alpha,\mathcal U]'=\phi(\Delta_{\alpha,\mathcal
U}\wedge\ker\pi)$ or, in other words, let
$[\alpha,\mathcal U]'$ be the unique compatible uniformity
on
$A$ such that $(\pi')^{-1}([\alpha,\mathcal
U]')=(\Delta_{\alpha,\mathcal
U}\wedge\ker\pi)\vee\ker\pi'$.

\begin{theorem}\label{T:Same} Let $A$ be an algebra, and
let
$\alpha\in\Con A$ and
$\mathcal U\in\OpUnif A$. Then
$[\alpha,\mathcal U]\leq[\alpha,\mathcal U]'$,
with equality if $A$ has an underlying group structure.
\end{theorem}

\begin{proof}
To show that $[\alpha,\mathcal
U]\leq[\alpha,\mathcal U]'$, it suffices to prove that
$\tilde C(\alpha,\mathcal U;[\alpha,\mathcal U]')$. We
have 
\begin{align*}
(\pi')^{-1}[\alpha,\mathcal
U]'&=(\Delta_{\alpha,\mathcal U}\wedge\ker\pi)\vee\ker\pi'
\\
&\geq\ker\pi'\circ(\Delta_{\alpha,\mathcal
U}\wedge\ker\pi)\circ\ker\pi'\\
&\geq\ker\pi'\circ(\mathcal
M(\alpha,\mathcal U)\wedge\ker\pi)\circ\ker\pi',
\end{align*}
which implies that
for all
$W\in[\alpha,\mathcal U]'$, there is a $Q\in\mathcal
M(\alpha,\mathcal U)$ such that
$\left(\begin{array}{cc}u_{11}&u_{12}\\u_{21}&u_{22}
\end{array}\right)\in Q$ and $u_{11}=u_{12}$ imply
$u_{21}\mathrel Wu_{22}$.
Given in addition $n$, $n'$, and $t$, there is a
$U_{n,n',t}\in\mathcal U$ such that $\mathrm
M_{n,n',t}(\alpha,U_{n,n',t})\subseteq Q$, by the
definition of $\mathcal M(\alpha,\mathcal U)$. It follows
that
$\left(\begin{array}{cc} u_{11} & u_{12} \\ u_{21} &
u_{22}
\end{array}\right)\in\mathrm
M_{n,n',t}(\alpha,U_{n,n',t})$ and 
$u_{11}=u_{12}$ imply
$u_{21}\mathrel Wu_{22}$, proving
that $\tilde C(\alpha,\mathcal
U;[\alpha,\mathcal U]')$.

Now assume $A$ has an underlying group structure. (In
particular, this implies that $A$ belongs to a
congruence-permutable variety and $\OpUnif A(\alpha)$ is
modular, by \cite[Theorem~6.2]{R02}.) To show
$[\alpha,\mathcal U]'\leq[\alpha,\mathcal U]$,
it suffices to show that
$\Delta_{\alpha,\mathcal U}\wedge\ker\pi\leq
(\pi')^{-1}([\alpha,\mathcal U])$.
For, this implies
\begin{align*}
(\pi')^{-1}([\alpha,\mathcal U]')
&=(\Delta_{\alpha,\mathcal U}\wedge\ker\pi)
\vee\ker\pi' \\
&\leq((\pi')^{-1}([\alpha,\mathcal
U])\wedge\ker\pi)
\vee\ker\pi'\\
&=(\pi')^{-1}([\alpha,\mathcal U]).
\end{align*}

To show $\Delta_{\alpha,\mathcal
U}\wedge\ker\pi\leq(\pi')^{-1}([\alpha,\mathcal U])$, it
suffices to show that
$\mathrm x_{\mathbf m}(\Delta_{\alpha,\mathcal
U})\leq[\alpha,\mathcal U]$. For, in that case, by
\cite[Lemma~7.1]{R02}, if
$W\in[\alpha,\mathcal U]$, and $\hat
W\in[\alpha,\mathcal U]$ is such that $\hat
W^{-1}\subseteq W$, then
$\bar W\in[\alpha,\mathcal U]$ can be chosen, and
$V\in\Delta_{\alpha,\mathcal U}$, such that $\mathrm
x_{\mathbf m}(V)\subseteq\bar W$, so that
$\left(\begin{array}{cc}y&z\\x&z\end{array}\right)\in
V\implies y\mathrel{\hat W}x\implies x\mathrel Wy$. We
have
$\ker\pi\wedge\ker\pi'=\bot\leq\Delta_{\alpha,\mathcal
U}$, and if
$\left(\begin{array}{cc}b&b\\x&y\end{array}\right)\in
\Delta_{\alpha,\mathcal U}$, we also have
$\left(\begin{array}{c}y\\x\end{array}\right)
\mathrel{\ker\pi'}
\left(\begin{array}{c}b\\x\end{array}\right)
\mathrel{\ker\pi}
\left(\begin{array}{c}b\\y\end{array}\right)
\mathrel{\ker\pi'}
\left(\begin{array}{c}y\\y\end{array}\right)
\mathrel{\ker\pi}
\left(\begin{array}{c}y\\x\end{array}\right)$.
It follows by the Shifting Lemma (\cite[Lemma~7.4]{R02})
that there is a $\bar V \in\Delta_{\alpha,\mathcal U}$
such that if
$\left(\begin{array}{cc}b&b\\x&y\end{array}\right)\in\bar
V$, then
$\left(\begin{array}{cc}y&y\\y&y\end{array}\right)\in V$,
implying $x\mathrel Wy$. It follows that $\bar
V\cap\ker\pi\subseteq(\pi')^{-1}(W)\cap\ker\pi$, proving
that
$\Delta_{\alpha,\mathcal
U}\wedge\ker\pi\leq(\pi')^{-1}([\alpha,\mathcal U])$.

Now, $\mathrm x_{\mathbf m}(\mathcal M(\alpha,\mathcal
U))=\mathcal X_{\mathbf m}(\alpha,\mathcal
U)\leq[\alpha,\mathcal U]$ by Theorem~\ref{T:Commutator}.
For ordinal numbers $\nu$, we define inductively
\[\mathcal R_\nu=\begin{cases}
\mathcal M(\alpha,\mathcal U), & \nu=0\\
\mathcal R_{\nu'}\circ\mathcal R_{\nu'}, & \nu=\nu'+1\\
\bigcap_{\nu'<\nu}\mathcal R_{\nu'}, &\nu \text{ a limit
ordinal.}
\end{cases}
\]

We claim that $\mathrm x_{\mathbf m}(\mathcal
R_\nu)\leq[\alpha,\mathcal U]$ for every ordinal $\nu$.
It is easy to see that $\mathcal R_\nu$ is a singly
compatible semiuniformity for every $\nu$, and that the
sequence becomes stationary at $\Delta_{\alpha,\mathcal
U}=\Ug\mathcal M(\alpha,\mathcal U)$. If there is a first
ordinal $\nu$ such that
$\mathrm x_{\mathbf m}(\mathcal
R_\nu)\not\leq[\alpha,\mathcal U]$, then clearly, $\nu\neq
0$ and $\nu$ is not a limit ordinal. Thus, to prove the
claim, and that $\mathrm x_{\mathbf
m}(\Delta_{\alpha,\mathcal U})\leq[\alpha,\,\mathcal U]$,
it suffices to show that if
$\mathrm x_{\mathbf m}(\mathcal R_\nu)\leq[\alpha,\mathcal
U]$, then
$\mathrm x_{\mathbf m}(\mathcal
R_{\nu+1})\leq[\alpha,\mathcal U]$.
We must show that, given $W\in[\alpha,\mathcal U]$, there
is an $R\in\mathcal R_{\nu+1}$ such that $\mathrm
x_{\mathbf m}(R)\subseteq W$.

We use the fact that the terms
$m_0(x,y,z,w)=w$, $m_1(x,y,z,w)=xz^{-1}yx^{-1}w$, and
$m_2(x,y,z,w)=w$ are a sequence of Day terms for any
variety of algebras with group structures.
 Without loss of generality,
by the discussion of Section~\ref{S:Groups}, we can
also assume that
$W$ is left translation invariant. Let
$\bar W\in[\alpha,\mathcal U]$ be left translation
invariant and such that $a\mathrel{\bar W}a'$ and
$b\mathrel{\bar W}b'\implies ab\mathrel Wa'b'$.
Since $\mathrm x_{\mathbf m}(\mathcal
R_\nu)\leq[\alpha,\mathcal U]$ by the induction
hypothesis, there is an $\bar R\in\mathcal R_\nu$ such
that $\mathrm x_{\mathbf m}(\bar R)\subseteq\bar W$. We
have
\[
\mathrm x_{\mathbf m}(\bar R\circ\bar R)=\left\{\,\pair
{xw^{-1}zx^{-1}u}x\mathrel{\bigg|}
\left(\begin{array}{cc}x&y\\u&v\end{array}\right),
\left(\begin{array}{cc}y&z\\v&w\end{array}\right)\in
\bar R\,\right\}\cup\left\{\,\pair xx\mathrel{\bigg|}
x\in A\,\right\},
\]
where the second term in the union takes care of the
contributions to $\mathrm x_{\mathbf m}(\mathcal
R_{\nu+1})$ coming from the terms $m_0$ and $m_2$. Since
$\mathrm x_{\mathbf m}(\bar R)\subseteq\bar W$, we have
$\pair{xv^{-1}yx^{-1}u}x$,
$\pair{yw^{-1}zy^{-1}v}y\in\bar W$. This implies that
$\pair{v^{-1}yx^{-1}u}e$, $\pair{w^{-1}zy^{-1}v}e\in\bar
W$, which implies that
\[
\pair{(w^{-1}zy^{-1}v)(v^{-1}yx^{-1}u)}e=\pair
{w^{-1}zx^{-1}u}e\in W.
\]
But $W$ was assumed left
translation invariant, so this implies that
$\pair{xw^{-1}zx^{-1}v}x\in W$, proving that $\mathrm
x_{\mathbf m}(\bar R\circ\bar R)\subseteq W$, and by
induction that $x_{\mathbf m}(\Delta_{\alpha,\mathcal
U})\leq[\alpha,\mathcal U]$.
\end{proof}

\section{Properties of the Commutator}\label{S:Properties}

In this section, we discuss general properties of the
commutator on uniformities, for a congruence-modular
variety.

\subsection{Elementary properties}

\begin{theorem}\label{T:Elem} Let $A$ be an
algebra in a congruence-modular variety $\mathbf V$. We have
\begin{enumerate}
\item[(1)] If $\mathcal U$, $\mathcal V\in\OpUnif A$, then
$[\mathcal U,\mathcal V]=[\mathcal V,\mathcal U]$.

\item[(2)] The commutator is monotone, i.e., if
$\mathcal U\leq\mathcal U'$, then $[\mathcal U,\mathcal
V]\leq[\mathcal U',\mathcal V]$.

\item[(3)] $[\mathcal U,\mathcal V]\leq\mathcal
U\wedge\mathcal V$.

\item[(4)] If $B\in\mathcal V$,
$f:B\to A$ is a homomorphism, and $\mathcal U$, $\mathcal
V\in\OpUnif A$, then $[f^{-1}(\mathcal U),f^{-1}(\mathcal
V)]\leq f^{-1}([\mathcal U,\mathcal V])$.

\end{enumerate}
\end{theorem}

\begin{proof}
(1): This is just a restatement of part of
Theorem~\ref{T:Commutator}.

(2):  If $\mathcal W$ satisfies the $\mathcal U',\mathcal
V$-term condition, then it also satisfies the $\mathcal
U,\mathcal V$-term condition. Or, use the obvious fact
that $\mathcal X_{\mathbf m}(\mathcal U,\mathcal
V)\leq\mathcal X_{\mathbf m}(\mathcal U',\mathcal V)$.

(3): Follows from Theorem~\ref{T:AllComm}.

(4): We have
\[\mathcal X_{\mathbf m}(f^{-1}(\mathcal
U),f^{-1}(\mathcal V))\leq f^{-1}(\mathcal X_{\mathbf
m}(\mathcal U,\mathcal V))\]
because $f$ is a homomorphism.
\end{proof}

\subsection{Additivity}

\begin{theorem}\label{T:Additive} Let $A$ be an algebra
in a congruence-modular variety with sequence of Day terms
$\mathbf m$, and let $\mathcal U_i$, for $i\in I$, and
$\mathcal V$ belong to $\OpUnif A$. Then
$\bigvee_{i\in I}[\mathcal U_i,\mathcal V]\leq
[\bigvee_{i\in I}\mathcal U_i,\mathcal V]$, with
equality if
$\mathcal V$ is of the form $\Ug\{\,\alpha\,\}$ and $A$
has an underlying group structure, or if
$I$ is finite and the
$\mathcal U_i$ permute pairwise.
\end{theorem}

\begin{proof} $\bigvee_i[\mathcal U_i,\mathcal
V]\leq[\bigvee_i\mathcal U_i,\mathcal V]$ by
Theorem~\ref{T:Elem}(2).

Let us denote
$\bigvee_i[\mathcal U_i,\mathcal V]$ by $\mathcal W$
for the remainder of the proof.

Suppose $\mathcal V=\Ug\{\,\alpha\,\}$ for some
congruence $\alpha$ and $A$ has an underlying group
structure. The commutator is commutative so we can switch
the arguments on each side. Since
$\alpha$ centralizes $\mathcal
U_i$ modulo
$[\alpha,\bigvee_i\mathcal U_i]$, we have
$\Delta_{\alpha,\mathcal U_i}\wedge\ker\pi\leq\mathcal
Z$, where $\mathcal
Z=(\pi')^{-1}([\alpha,\bigvee_i\mathcal U_i])$. Because
compatible extension preserves joins, we have
$\Delta_{\alpha,\vee_i\mathcal
U_i}=\bigvee_i(\Delta_{\alpha,\mathcal
U})=\Ug(\bigcap_i\Delta_{\alpha,\mathcal U_i})$, where the
singly compatible semiuniformity
$\mathcal R=\bigcap_i\Delta_{\alpha,\mathcal U_i}$
satisfies $\mathcal R\wedge\ker\pi\leq\mathcal Z$. Then
by the same argument used in the proof of
Theorem~\ref{T:Same},
$\Ug(\mathcal R)\wedge\ker\pi\leq\mathcal Z$. It follows
that
$[\alpha,\bigvee_i\mathcal U_i]\leq\mathcal W$.

Now suppose that $I=\{\,1,\ldots,k\,\}$ and that the
$\mathcal U_i$ permute pairwise. We will show that
$\mathcal W$ satisfies the
$(\bigvee_i\mathcal U_i),\mathcal V$-term condition.
Given $W\in\mathcal W$ and $t$, we define $W_0=W$ and for
each $i\in I$, inductively, define $U_i\in\mathcal U_i$,
$V_i\in\mathcal V$, and $W_i\in\mathcal W$ to be such
that $a\mathrel{U_i}a'$, $\mathbf b\mathrel{V_i}\mathbf
c$, and $t(a,\mathbf b)\mathrel{W_i}t(a,\mathbf c)$ imply
$t(a',\mathbf b)\mathrel{W_{i-1}}t(a',\mathbf c)$, using
the fact that $C(\mathcal U_i,\mathcal V;\mathcal W)$.
Then $a\mathrel{U_k\circ\ldots\circ U_1}a'$, $\mathbf
b\mathrel{(\bigcap_iV_i)}\mathbf c$, and $t(a,\mathbf
b)\mathrel{W_k}t(a,\mathbf c)$ imply $t(a',\mathbf
b)\mathrel Wt(a',\mathbf c)$. However,
$U_k\circ\ldots\circ U_1\in\bigvee_i\mathcal U_i$,
because the $\mathcal U_i$ permute pairwise; it
follows that $\mathcal W$ satisfies the
$(\bigvee_i\mathcal U_i),\mathcal V$-term condition,
proving that $\mathcal W=\bigvee_i[\mathcal U_i,\mathcal
V]=[\bigvee_i\mathcal U_i,\mathcal V]$.
\end{proof}

\begin{corollary}
The commutator of compatible uniformities on an algebra in
a congruence-permutable variety is finitely additive.
\end{corollary}

\begin{proof}
By \cite[Theorem~6.4]{R02},
compatible uniformities on an algebra in a
con\-gru\-ence-permut\-a\-ble variety permute pairwise.
\end{proof}

\subsection{The homomorphism property}

Recall the definition of $R(\mathcal W,U)$ from
Section~\ref{S:Compat}.

\begin{lemma}\label{T:X} Let $\mathbf V$ be a
congruence-modular variety with sequence of Day terms
$\mathbf m$. 
 If $A\in\mathbf V$, $\mathcal W\in\OpUnif A$,
$U$ and
$V$ are relations on $A$, and $n$, $n'$, and $t$ are
given, then
\[\mathrm x_{\mathbf m}(\mathrm M_{n,n',t}(R(\mathcal
W,U),R(\mathcal W,V)))\subseteq R(\mathcal W,\mathrm
x_{\mathbf m}(\mathrm M_{n,n',t}(U,V))).
\]
\end{lemma}

\begin{proof}
An element of $R(\mathcal
W,U)$ is a pair of equivalence classes of Cauchy nets with
respect to
$\mathcal W$, having
representatives such that for large enough indices, the
values taken by the representatives are related by
$U$. $R(\mathcal W,V)$ is defined similarly. Then an
element of $\mathrm x_{\mathbf m}(M_{n,n',t}(R(\mathcal
W,U),R(\mathcal W,V)))$ is a pair of equivalence classes
having representatives such that for large enough
indices, the values define pairs in
$\mathrm x_{\mathbf m}(M_{n,n',t}(U,V))$. That is, such a
pair belongs to
$R(\mathcal W,\mathrm x_{\mathbf m}((M_{n,n',t}(U,V)))$.
\end{proof}

\begin{theorem}
Let $\mathcal U$, $\mathcal V$, $\mathcal
W\in\OpUnif A$, where $A\in\mathbf V$, a
congruence-modular variety with Day terms $\mathbf m$.
Then
\begin{enumerate}
\item[(1)] If $\mathcal U\geq\mathcal W$ and $\mathcal
V\geq\mathcal W$, $[\mathcal U,\mathcal V]\vee\mathcal
W=(\nat\mathcal W)^{-1}[\mathcal U/\mathcal
W,\mathcal V/\mathcal W]$.
\item[(2)] If $\mathcal U$ and $\mathcal V$ permute with
$\mathcal W$, then $[\mathcal U,\mathcal V]\vee\mathcal
W=(\nat\mathcal W)^{-1}[(\mathcal U\vee\mathcal
W)/\mathcal W,(\mathcal V\vee\mathcal W)/\mathcal W]$.
\end{enumerate}
\end{theorem}

\begin{proof}
(2) follows from (1) because if $\mathcal U$ and $\mathcal
V$ permute with $\mathcal W$, then $[\mathcal U,\mathcal
V]\vee\mathcal W=[\mathcal U\vee\mathcal W,\mathcal
V\vee\mathcal W]\vee\mathcal W$, and we can apply (1)
since
$\mathcal U\vee\mathcal W\geq\mathcal W$ and $\mathcal
V\vee\mathcal W\geq\mathcal W$.

As for (1), first we have by Theorem~\ref{T:Elem}(4),
\begin{align*}
[\mathcal U,\mathcal V]&=[(\nat\mathcal
W)^{-1}(\mathcal U/\mathcal W),(\nat\mathcal
W)^{-1}(\mathcal V/\mathcal W)]\\
&\leq(\nat\mathcal W)^{-1}[\mathcal U/\mathcal
W,\mathcal V/\mathcal W]
\end{align*}
and also $\mathcal W=(\nat\mathcal W)^{-1}(\mathcal
W/\mathcal W)\leq(\nat\mathcal W)^{-1}[\mathcal
U/\mathcal W,\mathcal V/\mathcal W]$, because
$\mathcal W/\mathcal W$ is the least element of the
lattice $\OpUnif A/\mathcal W$, whence
\[
[\mathcal
U,\mathcal V]\vee\mathcal W\leq(\nat\mathcal
W)^{-1}([\mathcal U/\mathcal W,\mathcal V/\mathcal W]).
\]

To prove the opposite inequality, it suffices to show
that $\mathcal X_{\mathbf m}(\mathcal U/\mathcal
W,\mathcal V/\mathcal W)\leq([\mathcal U,\mathcal
V]\vee\mathcal W)/\mathcal W$.
If $Q\in[\mathcal U,\mathcal V]\vee\mathcal W$, then
in particular, $Q\in[\mathcal U,\mathcal V]$. Then by the
definition of $[\mathcal U,\mathcal V]$ and $\mathcal
X_{\mathbf m}(\mathcal U,\mathcal V)$, given $n$, $n'$,
and $t$, there exist $U\in\mathcal U$, $V\in\mathcal V$
such that $\mathcal X_{\mathbf m,n,n',t}(U,V)\subseteq Q$.
Then by Lemma~\ref{T:X}, $\mathcal X_{\mathbf
m,n,n',t}(R(\mathcal W,U),R(\mathcal W,V))\subseteq
R(\mathcal W,Q)$.
 This proves
that
$\mathcal X_{\mathbf m}(\mathcal U/\mathcal W,\mathcal
V/\mathcal W)\leq([\mathcal U,\mathcal V]\vee\mathcal
W)/\mathcal W$, because $R(\mathcal
W,U)\in\mathcal U/\mathcal W$ and $R(\mathcal
W,V)\in\mathcal V/\mathcal W$, and relations of the form
$R(\mathcal W,Q)$ for
$Q\in[\mathcal U,\mathcal V]\vee\mathcal W$ form a base
for $([\mathcal U,\mathcal V]\vee\mathcal W)/\mathcal W$.
\end{proof}

\section{Commutator Operations and Congruential
Uniformities}\label{S:Congruential}

Recall that a uniformity $\mathcal U\in\OpUnif A$ is
\emph{congruential} if it has a base of congruences.
Given a filter of congruences $F$, it determines a
congruential uniformity $\Ug F$, of which $F$ is
a base and which determines $F$.

For example, consider filters in $\Con\Z$, where
$\Z$ is the ring of integers.  In addition to the
principal filters $\Fg\{\,(\,n\,)\,\}$ for $n\in\N$, where
$\N$ is the set of natural numbers, there are many other
filters such as
$(\,m^\infty\,)\overset{\text{Def}}=\Fg\{\,(\,p^n\,)\mid
n\in\N\,\}$ for $m$ a nonzero natural number. (Of course,
the mapping $m\mapsto(m^\infty)$ is not one-one.)

The mapping $F\mapsto\Ug F$, from the lattice of
congruential uniformities of $A$ into the lattice of
uniformities, preserves arbitrary meets, and by
\cite[Theorem 6.3]{R02}, if the algebra $A$ has permuting
uniformities, it preserves finite joins.

For the time being, we will assume that the algebras we
are discussing belong to a congruence-modular variety.

\begin{proposition}\label{T:CUBase}
Let $\alpha\in\Con A$, and
$F$,
$F'\in\Fil\Con A$. Then the sets
$\{\,[\alpha,\phi]\mid\phi\in
F\,\}$,
$\{\,[\phi,\alpha]\mid\phi\in
F\,\}$, and
$\{\,[\phi,\phi']\mid\phi\in F,\,\phi'\in F'\,\}$
are bases for filters $[\alpha,F]$, $[F,\alpha]$,
and $[F,F']$ in $\Con A$.
\end{proposition}

\begin{proposition}\label{T:FCommF}
Let $\alpha$, $\beta\in\Con A$, and
$F\in\Fil\Con A$. Then
\begin{enumerate}
\item[(1)]
$[\Fg^{\Con A}\{\,\alpha\,\},\Fg^{\Con A}\{\,\beta\,\}]=
\Fg^{\Con A}\{\,[\alpha,\beta]\,\}$
\item[(2)]
$[\alpha,F]=[\Fg^{\Con A}\{\,\alpha\,\},F]$
\item[(3)]
$[F,\alpha]=[F,\Fg^{\Con A}\{\,\alpha\,\}]$.
\end{enumerate}
\end{proposition}

Thus, we have defined a binary operation on $\Fil\Con A$
which extends the commutator on $\Con A$. Clearly this
operation satisfies the elementery properties of the
commutator as given in Theorem~\ref{T:Elem}; we leave the
statement of the theorem to the reader. Let us also prove
that this commutator on $\Fil\Con A$ is finitely additive:

\begin{theorem}\label{T:CongAdditive} If $F_i$, $i=1$,
$2$,
$\ldots$,
$n$ and
$F'\in\Fil\Con A$, then
$[\bigvee_iF_i,F']=\bigvee_i[F_i,F']$.
\end{theorem}

\begin{proof}
Clearly we have $\bigvee_i[F_i,F']\leq[\bigvee_iF_i,F']$
by monotonicity. To prove the opposite inequality, let
$\chi\in\bigvee_i[F_i,F']$ be of the form
$\chi=\bigvee_i[\alpha_i,\beta_i]$ for $\alpha_i\in F_i$
and $\beta_i\in F'$. Then by the monotonicity and
additivity of the commutator on congruences, we have
$\chi\geq[\bigvee_i\alpha_i,\bigwedge_i\beta_i]\in
[\bigvee_iF_i,F']$. This proves that
$[\bigvee_iF_i,F']\leq\bigvee_i[F_i,F']$.
\end{proof}

Now, let us relate this commutator operation on $\Fil\Con
A$ to that on $\OpUnif A$:

\begin{theorem}\label{T:CUCommU}
Let $A$ be an algebra. We have
\begin{enumerate}
\item[(1)] If $\alpha\in\Con A$, and $F\in\Fil\Con A$,
then
$[\Ug\{\,\alpha\,\},\Ug F]\leq
\Ug[\alpha,F]$
\item[(2)] If $\alpha\in\Con A$, and $F\in\Fil\Con A$,
then
$[\Ug F,\Ug\{\,\alpha\,\}]\leq
\Ug [F,\alpha]$.
\item[(3)] If $F$, $F'\in\Fil\Con A$, then $[\Ug
F,\Ug F']\leq\Ug [F,F']$.
\end{enumerate}
\end{theorem}

\begin{proof} 
(3): If $W\in\Ug [F,F']$, then
$[\alpha,\beta]\subseteq W$ for some $\alpha\in F$ and
$\beta\in F'$, because such congruences form a base of
$\Ug[F,F']$.  $[\alpha,\beta]$ satisfies the
$\alpha,\beta$-term condition, so for all $t$, $a$ and
$a'$ such that $a\mathrel\alpha a'$, and $\mathbf b$ and
$\mathbf c$ with $\mathbf b\mathrel\beta\mathbf c$,
$t(a,\mathbf b)\mathrel{[\alpha,\beta]}t(a,\mathbf c)$
implies $t(a',\mathbf
b)\mathrel{[\alpha,\beta]}t(a',\mathbf c)$. But
$\alpha\in\Ug F$ and $\beta\in\Ug F'$.
Thus, $\Ug [F,F']$
satisfies the $\Ug F,\Ug F'$-term
condition. If follows that $[\Ug F,\Ug
F']\leq\Ug [F,F']$.

(1): By (3),
$[\Ug\{\,\alpha\,\},\Ug F]=[\Ug\Fg^{\Con A}\{\,\alpha\,\},
\Ug F]\leq\Ug [\Fg^{\Con A}\{\,\alpha\,\},F]=\Ug
[\alpha,F]$.

(2): similar to proof of (1).
\end{proof}

\begin{remark} Propositions~\ref{T:CUBase} and
\ref{T:FCommF}, and Theorem~\ref{T:CUCommU} hold in
non-congruence-modular varieties, if we replace
$[\alpha,\beta]$ by $C(\alpha,\beta)$ and define
$C(\alpha,F)$, $C(F,\alpha)$, and $C(F,F')$ or similarly
if we replace $[\alpha,\beta]$ by $\tilde
C(\alpha,\beta)$ and define $\tilde C(\alpha,F)$, $\tilde
C(F,\alpha)$, and $\tilde C(F,F')$. We omit the details.
\end{remark}

\section{Miscellany}\label{S:Examp}

\subsection{Congruential uniformities on
commutative rings}
 For
$A$ a commutative ring, Theorem~\ref{T:CUCommU} can
be improved. For a
translation invariant relation
$U$ on $A$, we denote by
$\delta(U)$ the set of differences $a-b$ for $a$,
$b\in A$ such that
$a\mathrel Ub$.

\begin{theorem} If $A$ is a commutative ring, then the
inequalities in the conclusion of Theorem~\ref{T:CUCommU}
are equalities.
\end{theorem}

\begin{proof} It suffices to show that $\Ug [F,F']\leq
[\Ug F,\Ug F']$. We use the fact that $[\Ug F,\Ug F']$
satisfies the $\Ug F,\Ug F'$-term condition for the
term $t(x,y)=xy$. Let
$U\in [\Ug F,\Ug F']$. By Lemma~\ref{T:InvBase}, we
may assume $U$ is translation-invariant without loss of
generality. There are $\bar U\in [\Ug F,\Ug F']$,
$\alpha\in F$, and $\beta\in F'$ such that
$a\mathrel\alpha a'$, $b\mathrel\beta c$, and
$ab\mathrel{\bar U}ac$ imply $a'b\mathrel Ua'c$. Then if
$b\mathrel\beta c$ and $a\in I_\alpha$ (the ideal
corresponding to $\alpha$),
$0b=0c$, implying
$ab\mathrel U ac$. It follows that
$I_{[\alpha,\beta]}=I_\alpha I_\beta\subseteq\delta(U)$,
and this implies that
$[\alpha,\beta]\subseteq U$.
\end{proof}

It follows from this theorem that, for commutative
rings, the commutators $[\alpha,F]$, $[F,\alpha]$ and
$[F,F']$ can be
considered as commutators of
uniformities.

\begin{example}
In $\OpUnif\mathbb Z$, we have
$[\Ug(p^\infty),\Ug(q^\infty)]=\Ug((pq)^\infty)$, for
prime numbers $p\neq q$, showing that the commutator of
two compatible uniformities is not always equal to
$\Ug\{\,\alpha\,\}$ for some congruence $\alpha$.
\end{example}

\subsection{Algebras in congruence-distributive varieties}

\begin{theorem}\label{T:ConDist}
Let $A$ be an algebra in a congruence-distributive
variety, with J\'onsson terms $d_0$, $\ldots$, $d_k$, and
let $\mathcal U$, $\mathcal V\in\OpUnif A$. Then
$[\mathcal U,\mathcal V]=\mathcal U\wedge\mathcal V$.
\end{theorem}

\begin{proof} We already know that $[\mathcal U,\mathcal
V]\leq\mathcal U\wedge\mathcal V$. To prove the opposite
inequality, we must show that if $W\in[\mathcal
U,\mathcal V]$, then there are $U\in\mathcal U$,
$V\in\mathcal V$ such that $U\cap V\subseteq W$. We will
use the J\'onsson terms to prove this.

Let $W\in[\mathcal U,\mathcal V]$. We
define $W_k$, $W_{k-1}$, $\ldots$, $W_0\in[\mathcal
U,\mathcal V]$ successively, as follows:
Set $W_k=W$. If $n$ is odd, then there exist
$W_{n-1}\in[\mathcal U,\mathcal V]$, $U_{n-1}\in\mathcal
U$, $V_{n-1}\in\mathcal V$ such that
$d_n(a,b,a)\mathrel{W_{n-1}}d_n(a,b,b)$ and
$a\mathrel{U_{n-1}\cap V_{n-1}}b$ imply
$d_n(a,b,a)\mathrel{W_n}d_n(a,b,b)$. If $n>0$ is even,
then there exist $W_{n-1}$, $U_{n-1}$, $V_{n-1}$ such that
$d_n(a,a,a)\mathrel{W_{n-1}}d_n(a,a,b)$ and
$a\mathrel{U_{n-1}\cap V_{n-1}}b$ imply
$d_n(a,a,a)\mathrel{W_n}d_n(a,a,b)$. Now let
$U=\bigcap_nU_n$, $V=\bigcap_nV_n$, and $a\mathrel{U\cap
V}b$. We have $d_0(a,a,a)=a=d_0(a,a,b)$ so
$d_0(a,a,a)\mathrel{W_0}d_0(a,a,b)$. We further have
\begin{align*}
d_2(a,b,a)=a=d_1(a,b,a)&\mathrel{W_1}d_1(a,b,b)
=d_2(a,b,b),\\
d_3(a,a,a)=a=d_2(a,a,a)&\mathrel{W_2}d_2(a,a,b)
=d_3(a,a,b),
\end{align*}
and so on, ending with
$a=d_k(a,b,a)\mathrel{W_k}d_k(a,b,b)=b$ if $k$ is even
and with $a=d_k(a,a,a)\mathrel{W_k}d_k(a,a,b)=b$ if $k$
is odd. In either case, we have shown that $U\cap
V\subseteq W_k=W$.
\end{proof}

\subsection{Abelian algebras}

Recall that an algebra $A$ in a congruence-modular
variety is \emph{abelian} if $[\top_A,\top_A]=\bot_A$.
In this case, from Theorem~\ref{T:Elem} and
Theorem~\ref{T:CommCong}, we have
$[\mathcal U,\mathcal V]=\Ug\{\,\bot_A\,\}$ for any
$\mathcal U$,
$\mathcal V\in\OpUnif A$.

For example, abelian groups
are abelian algebras, so we might consider the group of
real numbers and the commutator $[\mathcal U,\mathcal
U]$, where $\mathcal U$ is the unique
compatible uniformity (compatible, that is, with respect
to the abelian group operations) that gives rise to the
usual topology on the group of real numbers.
$\mathcal U$ is noncongruential and abelian. This example
shows that noncongruential uniformities can have a
commutator that is congruential, indeed of the form
$\Ug\{\,\alpha\,\}$ for $\alpha$ a congruence.

\section{Conclusions}\label{S:Conclu}

In this final section, we will review some of the
questions still open regarding the commutator of
uniformities, and uniform universal algebra generally.

The most important question is the possible additivity and
even complete additivity of the commutator, as holds for
congruences and as we have proved for some special cases
in Theorem~\ref{T:Additive} and
Theorem~\ref{T:CongAdditive}. Many applications of
commutator theory rely on this. An obstacle here is the
difficulty of dealing with joins of compatible
uniformities. A more specific question, which might be
easier to settle, is complete additivity for compatible
uniformities of an algebra in a congruence-permutable
variety. We proved finite additivity in this case, using
the fact that compatible uniformities permute pairwise.

The proof that, in the case of an algebra $A$ with
underlying group, formation of commutators with a
uniformity of the form
$\Ug\{\,\alpha\,\}$ is completely additive, utilizing the
theory of $\OpUnif A(\alpha)$, suggests that an
appropriate definition for
$A(\mathcal U)$ may help settle the additivity question.

 A uniformity
$\mathcal U$ on an algebra
$A$ in a congruence-modular variety $\mathbf V$ can be
defined as \emph{abelian} if $[\mathcal U,\mathcal
U]=\Ug\{\,\bot_A\,\}=\Fg\{\,\Delta_A\,\}$. In the case of
a congruence
$\alpha$, abelianness leads to a structure of abelian
group object on the algebra $A(\alpha)$, viewed as an
algebra over $A$ (that is, as an object in the comma
category of algebras of $\mathbf V$ over $A$). The
abelian group operations on this abelian group
object can be obtained from any difference term.
The problem of generalizing this theory to uniformities
again depends upon the proper definition of
$A(\mathcal U)$.

In the theory of uniform universal algebra, an important
open question is the possible modularity of the
uniformity lattice of an algebra in a congruence-modular
variety. This has only been proved for algebras in
congruence-permutable varieties and not more generally,
although there is a partial result
\cite[Theorem~7.5]{R02}.

A similar question is the
possible distributivity of the uniformity lattice of an
algebra in a congruence-distributive variety. This has
been proved only for arithmetic algebras
\cite[Theorem~6.5]{R02}. Note that because $[\mathcal
U,\mathcal V]=\mathcal U\wedge\mathcal V$ for congruence
distributive agebras (Theorem~\ref{T:ConDist}) additivity
of the commutator of compatible uniformities in  this
case is equivalent to distributivity of $\OpUnif A$, and
complete additivity is equivalent to the
distributivity of meet over an arbitrary join.

\subsection*{Acknowledgement}
We would like to thank Keith Kearnes for his ideas on
this topic, and some lively discussions.


\begin{thebibliography}{99}

\bibitem{B}
Nicholas Bourbaki, \emph{General Topology, Chapters
1--4.}  English translation of \emph{Topologie
G\'en\'eral}.
Springer--Verlag, Berlin, Heidelberg, New York, 1989.

\bibitem{B-S}
Stanley Burris and H. P. Sankappanavar,
\emph{A Course in Universal Algebra},
Graduate Texts in Mathematics 78,
Springer--Verlag, New York, 1981.

\bibitem{Day}
Alan Day,
\emph{A characterization of modularity for congruence
lattices of algebras},
Canad.\ Math.\ Bull.\ \textbf{12}(1969)
167--173.

\bibitem{F-McK}
Ralph Freese and Ralph McKenzie,
\emph{Commutator Theory
for
Congruence-Modular Varieties},
London
Math.\ Soc.\ Lecture Note Series 125,
Cambridge University
Press, 1987.

\bibitem{Jonsson}
Bjarni J\'onsson, \emph{Algebras whose congruence lattices
are distributive}, Math.\ Scand.\ \textbf{21}(1967)
110--121.

\bibitem{MacL}
Saunders Mac Lane,
\emph{Categories for the Working Mathematician},
Graduate Texts in Mathematics 5,
Springer--Verlag,
New York,
1971.

\bibitem{Mal}
A. I. Mal'tsev,
\emph{On the general theory of algebraic systems},
Mat.\ Sb.\ \textbf{35}(77), pp.\ 3--20.

\bibitem{ItoTopological}
\emph{Topological groups},
in \emph{Encyclopedic Dictionary of
Mathematics, Second Edition}, Kiyosi Ito, ed.,
The MIT Press,
Cambridge, Massachusetts-London,
1986.

\bibitem{ItoUniform}
\emph{Uniform spaces},
in \emph{Encyclopedic Dictionary of
Mathematics, Second Edition}, op.\ cit.

\bibitem{R02} William H.\ Rowan,
\emph{Algebras with a compatible
uniformity},
Algebra Universalis \textbf{47} (2002) 13--43.

\end{thebibliography}
\end{document}